\documentclass{article}
\pdfoutput=1
\usepackage{amssymb}
\usepackage{appendix}
\usepackage{latexsym}
\usepackage{mathrsfs}
\usepackage{enumerate}
\usepackage{extarrows}
\usepackage{graphicx}
\usepackage{caption}
\usepackage{supertabular}
\usepackage{subfigure}
\usepackage{setspace}
\usepackage{epstopdf}

\let\cal=\mathcal
\newtheorem{theo}{Theorem}[section]

\newtheorem{remark}[theo]{Remark}
\newtheorem{ques}{Question}
\newtheorem{lemma}[theo]{Lemma}
\newtheorem{claim}[theo]{Claim}

\newtheorem{con}[theo]{Conjecture}
\newtheorem{prop}[theo]{Proposition}
\newtheorem{fact}[theo]{Fact}
\newtheorem{defi}[theo]{Definition}

\def\q{\hspace*{\fill}$\Box$\medskip}

\def\endproofbox{\hskip 1.3em\hfill\rule{6pt}{6pt}}

{%
\noindent{\it Proof}
}%
{%
 \quad\hfill\endproofbox\vspace*{2ex}
}
\begin{document}
%\date{}
\title{An irrational Lagrangian density of a single hypergraph}
\author{Zilong Yan \thanks{School of Mathematics, Hunan University, Changsha 410082, P.R. China. Email: zilongyan@hnu.edu.cn.} \and Yuejian Peng \thanks{ Corresponding author. School of Mathematics, Hunan University, Changsha, 410082, P.R. China. Email: ypeng1@hnu.edu.cn. \ Supported in part by National Natural Science Foundation of China (No. 11931002).}
}

\maketitle
\begin{abstract}
The {\em Tur\'an number} of an $r$-uniform graph $F$, denoted by $ex(n,F)$, is the maximum number of edges in an $F$-free $r$-uniform graph on $n$ vertices. The {\em Tur\'{a}n density} of $F$ is defined as $\pi(F)=\underset{{n\rightarrow\infty}}{\lim}{ex(n,F) \over {n \choose r }}.$ Denote $\Pi_{\infty}^{(r)}=\{ \pi(\cal F): \cal F {\rm \ is \ a \ family \ of}\ r{\rm -uniform \ graphs}\},$ $\Pi_{fin}^{(r)}=\{ \pi(\cal F): \cal F {\rm \ is \ a \ finite \ family \ of}\ r{\rm -uniform \ graphs}\}$ and $\Pi_{t}^{(r)}=\{\pi(\cal F): \cal F {\rm \ is \ a \ family \ of } \ r{\rm-uniform \ graphs \ and} \\ \ |\cal F|\le t\}.$ For graphs, Erd\H{o}s-Stone-Simonovits (\cite{ESi}, \cite{ES}) showed that $\Pi_{\infty}^{(2)}=\Pi_{fin}^{(2)}=\Pi_{1}^{(2)}=\{0, {1 \over 2}, {2 \over 3}, ...,\\{l-1 \over l}, ...\}.$  We know quite few about the  Tur\'{a}n density of an $r$-uniform graph for $r\ge 3$.    Baber and Talbot \cite{BT}, and Pikhurko \cite{Pikhurko2} showed  that there is an irrational number in $\Pi_{3}^{(3)}$ and $\Pi_{fin}^{(3)}$ respectively, disproving a conjecture of Chung and Graham \cite{FG}. Baber and Talbot \cite{BT} asked whether $\Pi_{1}^{(r)}$ contains an irrational number. The Lagrangian of a hypergraph has been a useful tool in hypergraph extremal problems. The {\em Lagrangian density } of an $r$-uniform graph $F$ is $\pi_{\lambda}(F)=\sup \{r! \lambda(G):G\;is\;F\text{-}free\}$, where $\lambda(G)$ is the Lagrangian of an $r$-uniform graph $G$. Sidorenko \cite{Sidorenko-89} showed that  the Lagrangian density of an $r$-uniform hypergraph $F$ is the same as the  Tur\'{a}n density of the  extension of $F$.
 In this paper, we show that the Lagrangian density of $F=\{123, 124, 134, 234, 567\}$ (the disjoint union of $K_4^3$ and an edge)  is ${\sqrt 3\over 3}$, consequently, the Tur\'{a}n density of the extension of $F$ is an irrational number, answering  the question of Baber and Talbot.

\end{abstract}


\section{Introduction}
For a set $V$ and a positive integer $r$, let $V^{(r)}$ denote the family of all $r$-subsets of $V$. An {\em $r$-uniform graph} or {\em $r$-graph $G$} consists of a set $V(G)$ of vertices and a set $E(G) \subseteq V(G) ^{(r)}$ of edges. Let $e(G)$ denote the number of edges of $G$. An edge $e=\{a_1, a_2, \ldots, a_r\}$ will be simply denoted by $a_1a_2 \ldots a_r$. An $r$-graph $H$ is a {\it subgraph} of an $r$-graph $G$, denoted by $H\subseteq G$, if $V(H)\subseteq V(G)$ and $E(H)\subseteq E(G)$. A subgraph of $G$ {\em induced} by $V'\subseteq V$, denoted as $G[V']$, is the $r$-graph with vertex set $V'$ and edge set $E'=\{e\in E(G):e \subseteq V'\}$. For $S\subseteq V(G)$, let $G-S$ denote the subgraph of $G$ induced by $V(G)\setminus S$. Let $G^c$ denote the {\em complement} $r$-graph of an $r$-graph $G$ with $V(G^c)=V(G)$ and $E(G^c)=\{e: e\in V(G)^r\setminus E(G)\}$. Let $K^{r}_t$ denote the complete $r$-graph on $t$ vertices. Let $K^{r-}_t$ be obtained by removing one edge from $K_t^r$. For a positive integer $n$, let $[n]$ denote $\{1, 2, 3, \ldots, n\}$.

For a family ${\cal F}$ of $r$-graphs, an $r$-graph $G$ is called $\cal F$-free if it does not contain an isomorphic copy of any $r$-graph of $\cal F$. For a fixed positive integer $n$ and a family of $r$-graphs $\cal F$, the {\em Tur\'an number} of $\cal F$, denoted by $ex(n,\cal F)$, is the maximum number of edges in an $\cal F$-free $r$-graph on $n$ vertices. An averaging argument of Katona, Nemetz and Simonovits \cite{KNS} shows that the sequence ${ex(n,\cal F) \over {n \choose r }}$ is  non-increasing. Hence $\underset{n\rightarrow\infty}{\lim}{ex(n, \cal F) \over {n \choose r } }$ exists. The {\em Tur\'{a}n density} of $\cal F$ is defined as $$\pi(\cal F)=\lim_{n\rightarrow\infty}{ex(n, \cal F) \over {n \choose r }}.$$
If $\cal F$ consists of an single $r$-graph $F$, we simply write $ex(n, \{F\})$ and $\pi(\{F\})$ as $ex(n, F)$ and $\pi(F)$. Denote
$$\Pi_{\infty}^{(r)}=\{ \pi(\cal F): \cal F {\rm \ is \ a \ family \ of \ } r{\rm-uniform \ graphs} \}, $$
$$\Pi_{fin}^{(r)}=\{\pi(\cal F): \cal F {\rm \ is \ a \ finite \ family \ of \ } r{\rm-uniform \ graphs} \}$$
and
$$\Pi_{t}^{(r)}=\{ \pi(\cal F): \cal F {\rm \ is \ a \ family \ of \ } r{\rm-uniform \ graphs \ and }\ \vert \cal F \vert\le t  \}. $$
Clearly,
$$\Pi_{1}^{(r)}\subseteq \Pi_{2}^{(r)}\subseteq \cdots \subseteq\Pi_{fin}^{(r)}\subseteq \Pi_{\infty}^{(r)}.$$

For 2-graphs, Erd\H{o}s-Stone-Simonovits (\cite{ESi}, \cite{ES}) determined the Tur\'an numbers of all non-bipartite graphs asymptotically. Their result implies that $$\Pi_{\infty}^{(2)}=\Pi_{fin}^{(2)}=\Pi_{1}^{(2)}=\{0, {1 \over 2}, {2 \over 3}, ..., {l-1 \over l}, ...\}.$$
Very few results are known for $r\ge 3$. In \cite{FG} Chung and Graham  proposed the conjecture that every element in $\Pi_{fin}^{(r)}$ is a rational number. Baber and Talbot \cite{BT}, and Pikhurko \cite{Pikhurko2} disproved this conjecture by showing that there is an irrational number in $\Pi_{3}^{(3)}$ and $\Pi_{fin}^{(3)}$, respectively. %Baber and Talbot \cite{BT} applied Razborov's flag algebra framework to show that the Tur\'an density of a family of three $3$-graphs is the Lagrangian of a corresponding $3$-graph which is an irrational number. Independently, Pikhurko~\cite{Pikhurko2} showed that $\Pi_{fin}^{(r)}$ contains an irrational number for every $r\ge 3$ by applying the Strong Removal Lemma of R\"{o}dl and Schacht \cite{rodlsch}.
Baber and Talbot \cite{BT} asked whether $\Pi_{1}^{(r)}$ contains an irrational number. In this paper, we answer this question by showing that the Lagrangian density of the disjoint union of $K_4^3$ and an edge is an irrational number.

The hypergraph Lagrangian method has been helpful in  hypergraph extremal problems.
\begin{defi}
Let $G$ be an $r$-graph on $[n]$ and let
  $\vec{x}=(x_1,\ldots,x_n) \in [0,\infty)^n$. Define the {\em Lagrangian} function
$$\lambda (G,\vec{x})=\sum_{e \in E(G)}\prod\limits_{i\in e}x_{i}.$$
\end{defi}
The {\em Lagrangian} of
$G$, denoted by $\lambda (G)$, is defined as
 $$\lambda (G) = \max \{\lambda (G, \vec{x}): \vec{x} \in \Delta \},$$
where $$\Delta=\{\vec{x}=(x_1,x_2,\ldots ,x_n) : \sum_{i=1}^{n} x_i =1, x_i\geq0 \;for\;every\;i\in[n] \}.$$

The value $x_i$ is called the {\em weight} of the vertex $i$ and a vector $\vec{x} \in {\Delta}$ is called a {\em feasible weight vector} on $G$.
A feasible weight vector  $\vec{y}\in {\Delta}$ is called an {\em optimum weight vector} for $G$ if $\lambda (G, \vec{y})=\lambda(G)$.

In \cite{MS}, Motzkin and Straus established a connection between the Lagrangian of a $2$-graph and it's maximum complete subgraphs.
\begin{theo} {\em(\cite{MS})} \label{MStheo}
If $G$ is a $2$-graph in which a maximum complete subgraph has  $t$ vertices, then
$\lambda(G)=\lambda(K_t^2)={1 \over 2}(1 - {1 \over t})$.
\end{theo}

They also applied this connection to give another proof of the theorem of Tur\'an  on the Tur\'an density of complete graphs. Since then the Lagrangian method has been a useful tool in hypergraph extremal problems. Earlier applications include that Frankl and R\"{o}dl \cite{FR} applied it in disproving the long standing jumping constant conjecture of Erd\H{o}s. Sidorenko \cite{Sidorenko-89} applied Lagrangians of hypergraphs to first find infinitely many Tur\'an densities of  hypergraphs. More recent developments of this method were obtained in  \cite{Pikhurko, BT,  HK, NY, BIJ, NY2, Jenssen, JPW, HPW, YP}.  Determining the Lagrangian of a hypergraph is much more difficult than graphs and there is no conclusion similar to Theorem \ref{MStheo} for hypergraphs.
It is of great interests to estimate  Lagrangians of hypergraphs that have some certain properties. In 1980's, Frankl and F\"uredi \cite{FF} asked the question that for a given integer $m$,  what is  the maximum  Lagrangian among all $r$-graphs with $m$ edges?
  They conjectured that  the $r$-graph with $m$ edges formed by taking the first $m$ sets in the colex ordering of ${\mathbb N}^{(r)}$ has the largest Lagrangian of all $r$-graphs with  $m$ edges.  For distinct $A, B \in {\mathbb N}^{r}$ we say that $A$ is less than $B$ in the {\em colex ordering} if $max(A \triangle B) \in B$, where $A \triangle B=(A \setminus B)\cup (B \setminus A)$. By Theorem \ref{MStheo}, this conjecture is true when $r=2$. For hypergraphs,
Talbot \cite{T} first proved the conjecture for
 $r=3$ and ${\ell\choose 3}\leq m\leq {\ell\choose 3}+{{\ell-1}\choose 2}-\ell$, where $\ell>0$ is an integer. Subsequent progress in this conjecture  were made in the papers of Tang, Peng, Zhang and Zhao \cite{TPZZ1, TPZZ2}, Tyomkyn \cite{Tyo}, Lei, Lu and Peng \cite{LLP2018}, Nikiforov\cite{NIK}, Lei and Lu\cite{LL}, and Lu\cite{Lu}. Recently, Gruslys, Letzter and Morrison \cite{GLM2018} confirmed this conjecture for $r$ and ${\ell\choose r}\leq m\leq {\ell\choose r}+{{\ell-1}\choose r-1}$ if $\ell$ is sufficiently large. They also found infinitely many counterexamples for all $r\ge 4$. As remarked in \cite{GLM2018}, it would be  interesting to find the maximisers for other values of $m$ though it might be a very hard problem.
In this paper,  we will apply the connection of the Lagrangian density  and the  Tur\'an density of an $r$-graph to answer the question of Baber and Talbot. Our proof relies heavily on the estimation of  Lagrangians of $3$-graphs.

The {\em Lagrangian density } $\pi_{\lambda}(F)$ of  an $r$-graph $F$ is defined to be
$$\pi_{\lambda}(F)=\sup \{r! \lambda(G): G \:\: {\rm is} \:\: F{\text-}{\rm free}\}.$$
A pair of vertices $\{i, j\}$ is {\em covered}  in a hypergraph $F$ if there exists  an edge $e$ in $F$ such that  $\{i, j\}\subseteq e$. We say that $F$ covers pairs if every pair of vertices in $F$ is covered.
 Let $r\ge 3$ and $F$ be an $r$-graph.  The  {\em extension} of $F$, denoted by $H^F$ is
obtained as follows: For each pair of vertices $v_i$ and $v_j$ not covered in $F$, we add a set $B_{ij}$ of $r-2$ new vertices and the edge $\{v_i,v_j\} \cup B_{ij}$, where the $B_{ij}$'s are pairwisely disjoint over all such pairs $\{i,j\}$.

The Lagrangian density is closely related to the Tur\'an density.
The following proposition is implied by Theorem $2.6$ in \cite{Sidorenko-87} (see Proposition $5.6$ in  \cite{BIJ} and Corollary $1.8$ in \cite{Sidorenko-89} for the explicit statement). %It gives the  connection between the Lagrangian density of a hypergraph  $F$ and the Tur\'an density of its extension.
\begin{prop}{\em (\cite{Sidorenko-87, BIJ,Sidorenko-89})}\label{relationlt} Let $F$ be an $r$-graph. Then \\
$(i)$ $\pi(F)\le \pi_{\lambda}(F);$ \\
$(ii)$ $\pi(H^F)=\pi_{\lambda}(F).$ In particular, if $F$ covers pairs, then $\pi(F)= \pi_{\lambda}(F).$
\end{prop}

To answer the question of Baber and Talbot, we show that the Lagrangian density of $\{123, 124, 134, 234,\\567\}$,  the disjoint union of $K_4^3$ and an edge, denoted as $K_4^3\cup e$, is ${\sqrt 3\over 3}$. The following is our main Theorem.
\begin{theo}\label{main}
$\pi_{\lambda}(K_4^3\cup e)=\frac{\sqrt 3}{3}$.
\end{theo}
Applying Theorem \ref{main} and Proposition \ref{relationlt}, we see that the Tur\'an density of the extension of $K_4^3\cup e$ is ${\sqrt 3\over 3}$.

For an $r$-graph $H$ on $t$ vertices, it is clear that $\pi_{\lambda}(H)\ge r!\lambda{(K_{t-1}^r)}$.
An $r$-graph $H$ on $t$ vertices is {\em $\lambda$-perfect} if $\pi_{\lambda}(H)= r!\lambda{(K_{t-1}^r)}$.
Theorem \ref{MStheo} implies that all $2$-graphs are $\lambda$-perfect. Theorem \ref{main}  indicates that $K_4^3\cup e$ is not $\lambda$-perfect. We can show however that $K_4^3\cup k\cdot e$, the disjoint union of $K_4^3$ and $k$ disjoint edges is $\lambda$-perfect for $k\geq 2$.

\begin{theo}\label{perfect}
$K_4^3\cup k\cdot e$ is $\lambda$-perfect for $k\geq 2$.
\end{theo}

In Section \ref{sec2}, we give a sketch of the proof of Theorem \ref{main}. In Section \ref{perfect1}, we will give the proof of Theorem \ref{perfect}. In Section \ref{sec3}, we give some preliminaries on KKT conditions for continuous optimization problems and properties of hypergraph Lagrangians. In Section \ref{sec4}, we prove the main Lemmas needed in the proof of Theorem \ref{main}.
%The proofs of these Lemmas highly rely on the estimation of Lagrangians of hypergraphs. Our proofs in Section \ref{sec4} are theoretical. We also used Lingo to run corresponding optimization problems. The optimum values obtained by Lingo are also consistent with the expected optimum valus in our theoretical proof.

\section{Sketch of the proof of Theorem \ref{main}}\label{sec2}
The following three 3-graphs are to be used throughout the paper.

\textbf{B(2, n-2):} the 3-graph with vertex set $[n]$ and edge set $E(B(2, n-2))=\{e\in {[n]\choose 3}: e\cap\{1, 2\}\not=\emptyset\}$, i.e., every edge in $B(2, n-2)$ contains vertex $1$ or $2$ or both. Note that $B(2, n-2)$ is $K_4^3\cup e$-free, we will show that it is an extremal $3$-graph for $K_4^3\cup e$ (in terms of Lagrangian density).

$\mathbf{X_i:}$ the 3-graph with vertex set $[2i+2]$ such that $\{1, 2, 2j+1, 2j+2\}$ form $K_4^3$ for all $j$, $1\leq j\leq i$, i.e., it consists of $i$ copies of $K_4^3$ all sharing  vertices  $\{1, 2\}$.

$\mathbf{Y_i:}$ the 3-graph with vertex set $[i+3]$ such that $\{1, 2, 3, j+3\}$ form $K_4^3$ for all $j$, $1\leq j\leq i$, i.e., it consists of $i$ copies of $K_4^3$ all sharing  vertices  $\{1, 2, 3\}$.

An $r$-graph $G$ is {\em dense} if $\lambda (G') < \lambda (G)$ for every proper subgraph $G'$ of $G$.

{\em Sketch of the proof of Theorem \ref{main}:} For the lower bound, note that $B(2, n-2)$ is $K_4^3\cup e$-free, we shall prove $\underset{n\to\infty}{\lim}\lambda(B(2, n))={\sqrt 3\over 18}$ in Lemma \ref{B2n}. So $\pi_{\lambda}(K_4^3\cup e)\geq3!\underset{{n\rightarrow\infty}}{\lim}\lambda(B(2, n-2))={\sqrt3\over 3}$.

For the upper bound, let $G$ be a $K_4^3\cup e$-free 3-graph, our goal is to show that $\lambda(G)\leq\frac{\sqrt 3}{18}$. If $G$ is not dense, then there exists a proper subgraph $G'$ of $G$ such that $\lambda{(G')}=\lambda{(G)}$ and $|V(G')|<|V(G)|.$  If $G'$ is dense, then we stop. Otherwise, we continue this process until we find a dense subgraph $G''$ such that $\lambda(G'')=\lambda(G)$. This process terminates since the number of vertices is reduced by at least one in each step. To show $\lambda(G)\leq\frac{\sqrt 3}{18}$, it's sufficient to show that $\lambda(G'')\leq\frac{\sqrt 3}{18}$. So we may assume that $G$ is a dense $K_4^3\cup e$-free 3-graph. Suppose that $\lambda(G)>{\sqrt 3\over 18}$, we will prove the following lemmas in Section \ref{sec4}.

\begin{lemma}\label{X2}
Let $G$ be a dense $K_4^3\cup e$-free 3-graph with $\lambda(G)>{\sqrt 3\over 18}$. Then $G$ is $X_2$-free.
\end{lemma}

\begin{lemma}\label{2K4}
Let $G$ be a dense $K_4^3\cup e$-free 3-graph with $\lambda(G)>{\sqrt 3\over 18}$. Then $G$ contains at least two copies of $K_4^3$.
\end{lemma}

\begin{lemma}\label{Y2}
Let $G$ be a dense $K_4^3\cup e$-free 3-graph with $\lambda(G)>{\sqrt 3\over 18}$. Then $G$ is $Y_2$-free.
\end{lemma}

By Lemma \ref{2K4}, $G$ contains two copies of $K_4^3$.  Since $G$ is $K_4^3\cup e$-free, these two copies of $K_4^3$ must have $2$ or $3$ vertices in common. So $G$ contains a copy of $X_2$ or $Y_2$, a contradiction to Lemmas \ref{X2} and \ref{Y2}.

To complete the proof of Theorem \ref{main}, what remains is to show those three main lemmas. They will be given in Section \ref{sec4}.

\section{The proof of Theorem \ref{perfect}}\label{perfect1}
In order to prove Theorem \ref{perfect}, we need some lemmas form \cite{YP2}. Let $S_{2, t}$ denote the 3-graph with vertex set $\{v_1, v_2, u_1, u_2,..., u_t\}$ and edge set $\{v_1v_2u_1, v_1v_2u_2, ..., v_1v_2u_t\}$. A result of Sidorenko in \cite{Sidorenko-89} implies that $S_{2, t}$ is $\lambda$-perfect.

\begin{theo} \label{HS1}(\cite{YP2})
If $H$ is $\lambda$-perfect, then $H\cup S_{2, t}$ is $\lambda$-perfect for any $t\geq1$.
\end{theo}

\begin{claim}\label{3.1}(\cite{YP2})
Let $G$ be a 3-graph with $\lambda(G)>\lambda(K_{k+1}^3)$ and let $\vec{x}$ be an optimal weight vector. Then for any $v\in V(G)$, its weight $x_v$ satisfies that $x_v<1-\frac{\sqrt{k(k-1)}}{k+1}$.
\end{claim}

\begin{claim}\label{3.3}(\cite{YP2})
Let $v$ be a vertex in a 3-graph $G$ and $x_v$ be the weight of v in an optimal weight vector $\vec x$ of $G$. If $G-\{v\}$ is $H$-free, then $\lambda(G)\leq\frac{\pi_{\lambda}(H)(1-x_v)^3}{6(1-3x_v)}$.
\end{claim}

\begin{remark} \label{remark33}(\cite{YP2})
$f(x)=\frac{(1-x)^3}{1-3x}$ is increasing in $(0, \frac{1}{3}).$
\end{remark}

\begin{defi}
For $v\in V(G)$, the link graph of $v$ in $G$, denote by $G_v$, is the graph on vertex set $V(G)$ and the edge set $\{e\setminus\{v\}:v\in e\in E(G)\}$. Let $\omega(G_v)$ be the number of vertices in a maximum complete subgraph in $G_v$.
\end{defi}

\begin{claim}\label{3.5}(\cite{YP2})
Let a 3-graph $G$ be $H\cup S_{2, t}$-free, where H is a 3-graph with s vertices. Let $v\in V(H)$. If $H\subseteq G-\{v\}$, then $\omega(G_{v})\leq s+t$.
\end{claim}

\begin{claim}\label{3.4}(\cite{YP2})
Let a 3-graph G be $H\cup S_{2, t}$-free, where H is a 3-graph with s vertices. Let $S_{2, s+t}=\{v_1v_2b_1, v_1v_2b_2, ..., v_1v_2b_{s+t} \}\subseteq G$. Then $G-\{v_1, v_2\}$ is H-free.
\end{claim}

\begin{claim}\label{3.6}(\cite{YP2})
Let a 3-graph G be $H\cup S_{2, t}$-free, where H is a 3-graph with s vertices. If $H\subseteq G-\{v_1\}$ and $H\nsubseteq G-\{v_1, v_2\}$, then $\omega((G-\{v_2\})_{v_1})\leq s+t-1.$
\end{claim}

\emph{Proof of Theorem \ref{perfect}.} By Theorem \ref{HS1}, it's sufficient to show that $K_4^3\cup 2\cdot e$ is $\lambda$-perfect. Note that $K_4^3\cup 2\cdot e$ has 10 vertices. It's sufficient to show that if $G$ is $K_4^3\cup 2\cdot e$-free dense 3-graph then  $\lambda(G)\leq\lambda(K_9^3)$. Suppose on the contrary that $\lambda(G)>\lambda(K_9^3)={28 \over 243}$. Let $\vec x$ be an optimal weight vector of $G$.

\emph{Case 1.} There exists $v\in V(G)$ with weight $x_v$ such that $G-\{v\}$ is $K_4^3\cup e$-free.

By Claim \ref{3.1}, $x_v<1-\frac{2\sqrt {14}}{9}$. By Theorem \ref{main} and Claim \ref{3.3}, $$\lambda(G)\leq\frac{\frac{\sqrt 3}{18}(1-x_v)^3}{1-3x_v}=f(x_v).\eqno {(3)}.$$
Since $f(x_v)$ is increasing in $[0, 1-\frac{2\sqrt {14}}{9}]$, then
\begin{eqnarray*}
\lambda(G)&\leq&f(1-\frac{2\sqrt {14}}{9})\\
&=&\frac{28\sqrt{42}}{729(3\sqrt{14}-9)}\\
&\leq&\frac{28}{243},
\end{eqnarray*}
a contradiction.

\emph{Case 2.} For any $v\in V(G)$, $K_4^3\cup e\subseteq G-\{v\}$.

Since $\lambda(G)>\lambda(K_9^3)$ and $S_{2, 8}$ is $\lambda$-perfect, then $S_{2, 8}=\{v_1v_2b_1, v_1v_2b_2, ..., v_1v_2b_8\}\subseteq G$. By Claim \ref{3.4}, $G-\{v_1, v_2\}$ is $K_4^3\cup e$-free. Applying Claim \ref{3.5} ($s=7, t=1$), we have $$\omega(G_{v_1})\leq 8\quad{\rm and}\quad\omega(G_{v_2})\leq 8.$$ Applying Claim \ref{3.6} ($s=7, t=1$), we have $\omega((G-\{v_2\})_{v_1})\leq 7$ and $\omega((G-\{v_1\})_{v_2})\leq 7$.

Assume the weight of $v_1$ and $v_2$ are $a_1$ and $a_2$ respectively, and $a_1+a_2=2a$. Since $G-\{v_1, v_2\}$ is $K_4^3\cup e$-free and by Theorem \ref{main}, the contribution of edges containing neither $v_1$ nor $v_2$ to $\lambda(G, \vec x)$ is at most $\frac{\sqrt 3}{18}(1-2a)^3$. Since $\omega((G-\{v_2\})_{v_1})\leq 7$ and $\omega((G-\{v_1\})_{v_2})\leq 7$, by Theorem \ref{MStheo}, the contribution of edges containing either $v_1$ or $v_2$ to $\lambda(G, \vec x)$ is at most $2\times\frac{1}{2}a(1-\frac{1}{7})(1-2a)^2$. The contribution of edges containing both $v_1$ and $v_2$ to $\lambda(G, \vec x)$ is at most $a^2(1-2a)$. Therefore
\begin{eqnarray*}
\lambda(G)&\leq&\frac{\sqrt 3}{18}(1-2a)^3+a^2(1-2a)+\frac{6}{7}a(1-2a)^2\\
&\leq&\frac{1}{10}(1-2a)^3+a^2(1-2a)+\frac{6}{7}a(1-2a)^2=f(a)\\
f'(a)&=&\frac{66a^2-86a+9}{35}.
\end{eqnarray*}
Since $f(a)$ is increasing in $[0, \frac{43-\sqrt{1255}}{66}]$ and is decreasing in $[\frac{43-\sqrt{1255}}{66}, 1]$, then $\lambda(G)\leq f(\frac{43-\sqrt{1255}}{66})<\frac{28}{243}$.
\q

\section{Preliminaries}\label{sec3}
\subsection{Karush-Kuhn-Tucker Conditions}
Let us consider the optimisation problem:
\begin{flushleft}
\quad\quad\quad\quad\quad\quad\quad\quad maximise\ $f(x)$\\
\quad\quad\quad\quad\quad\quad\quad\quad subject to $g_i(x)\leq 0$, $i=1,\dots,m,$\quad\quad\quad\quad\quad\quad\quad\quad\quad\quad\quad\quad\quad\quad\quad\quad\quad\quad\quad(3.1)
\end{flushleft}
where $x\in \mathbb{R}^n$ and $f$ and $g_i$ are differentiable functions from $\mathbb{R}^n$ to $\mathbb{R}$ for all $i$. Let $\nabla{f(x)}$ be the gradient of $f$ at $x$ i.e. the vector in $\mathbb{R}^n$ whose $i$th coordinate is ${\partial\over\partial{x_i}}f(x)$. We say that KKT conditions hold at $x^*\in \mathbb{R}^n$ if there exist $\lambda_1,\dots\lambda_m\in \mathbb{R}$ such that
\begin{enumerate}[(i)]
\item $\nabla{f(x^*)}=\sum_{i=1}^m\lambda_i\nabla{g_i(x^*)},$
\item $\lambda_i\geq 0, i=1,\dots,m,$
\item $\lambda_ig_i(x^*)=0, i=1,\dots,m.$
\end{enumerate}
We call the constraints linear if $g_1,\dots,g_m$ are all affine functions.
\begin{theo}\label{KKT}(\cite{BV},\cite{Jenssen})
If the constraints of (3.1) are linear, then any optimum point of (3.1) must satisfy the KKT conditions.
\end{theo}
\subsection{Properties of the Lagrangian function}
The following fact follows immediately from the definition of the Lagrangian.
\begin{fact}\label{mono}
Let $G_1$, $G_2$ be $r$-graphs and $G_1\subseteq G_2$. Then $\lambda (G_1) \le \lambda (G_2).$
\end{fact}

\begin{fact} {\em (\cite{FR})}\label{fact2}
Let $G$ be an $r$-graph on $[n]$. Let $\vec{x}=(x_1,x_2,\dots,x_n)$ be an optimum weight vector on  $G$. Then
$$ \frac{\partial \lambda (G, \vec{x})}{\partial x_i}=r\lambda(G)$$
for every $i \in [n]$ satisfying $x_i>0$.
\end{fact}

Given an $r$-graph $G$, and $i, j\in V(G),$ define $$L_G(j\setminus i)=\{e: i\notin e, e\cup\{j\}\in E(G)\:and\: e\cup\{i\}\notin E(G)\}.$$

\begin{fact}\label{symmetry}
Let $G$ be an $r$-graph on $[n]$. Let $\vec{x}=(x_1,x_2,\dots,x_n)$ be a feasible weight vector on $G$. Let $i,j\in [n]$, $i\neq j$ satisfying $L_G(i \setminus j)=L_G(j \setminus i)=\emptyset$. Let
$\vec{y}=(y_1,y_2,\dots,y_n)$ be defined by letting $y_\ell=x_\ell$ for every $\ell \in [n]\setminus \{i,j\}$ and $y_i=y_j={1 \over 2}(x_i+x_j)$.
Then $\lambda(G,\vec{y})\geq \lambda(G,\vec{x})$. Furthermore, if the pair $\{i,j\}$ is contained in an edge of $G$, $x_i>0$ for each $1\le i\le n$,  and $\lambda(G,\vec{y})=\lambda(G,\vec{x})$, then $x_i=x_j$.
\end{fact}
{\em  Proof of Fact \ref{symmetry}.}
Since $L_G(i \setminus j)=L_G(j \setminus i)=\emptyset$, then
$$\lambda(G,\vec{y})-\lambda(G,\vec{x})=\sum_{\{i,j\} \subseteq e \in G}\left({(x_i+x_j)^2 \over 4}-x_ix_j\right)\prod\limits_{k\in e\setminus \{i,j\}}x_k \ge 0.$$
If the pair $\{i,j\}$ is contained in an edge of $G$ and $x_i>0$ for each $1\le i\le n$, then the equality holds only if $x_i=x_j$.
\q

\begin{fact}\label{*}
Let $\vec{x}=(x_1,x_2,\dots,x_n)$ be an optimum vector for an $r$-graph $G$ on $[n]$. If $L_G(j \setminus i)=\emptyset$, then we may assume that $x_i\geq x_j$.
\end{fact}
{\em  Proof of Fact \ref{*}.} If $x_i<x_j$, then let $\epsilon={x_j-x_i\over 2}$ and $\vec{x'}=(x_1,x_2,\dots,x_i+\epsilon,\dots,x_j-\epsilon,\dots,x_n)$. Since $L_G(j \setminus i)=\emptyset$, then
$$\lambda(G,\vec{x'})-\lambda(G,\vec{x})\geq\sum_{\{i,j\} \subseteq e \in G}((x_i+\epsilon)(x_j-\epsilon)-x_ix_j)\prod\limits_{k\in e\setminus \{i,j\}}x_k \geq {(x_j-x_i)^2\over 4}\sum_{\{i,j\} \subseteq e \in G}\prod\limits_{k\in e\setminus \{i,j\}}x_k\geq0.$$
\q
\begin{fact} {\em (\cite{FR})}\label{dense}
Let $G=(V,E)$ be a dense $r$-graph. Then $G$ covers pairs.
\end{fact}
Let $C_{r, m}$ denote the $r$-graph with $m$ edges formed by taking the first $m$ sets in the colex ordering of $\mathbb{N}^{r}$. The following two results was given in \cite{TPP} and \cite{T}, respectively.
\begin{lemma}\label{TPP}(\cite{TPP})
Let $m$ and $t$ be positive integers satisfying the condition that ${t\choose 3}-6\leq m\leq {t\choose 3}-3$. Let $G$ be a 3-graph with $m$ edges. Then $\lambda(G)\leq\lambda(C_{3, m})$.
\end{lemma}
\begin{lemma}\label{T}(\cite{T})
For any integers $m$, $t$ and $r$ satisfying the condition that ${t-1\choose r}\leq m\leq {t-1\choose r}+{t-2\choose r-1}$, then we have $\lambda(C_{r, m})=\lambda(K_{t-1}^r)$.
\end{lemma}
\begin{fact}\label{aaa}
Let $f(x)=\frac{x^2(1-x)}{4}+\frac{x(1-x)^2}{2}$, where $0\leq x\leq 1$. Then $f(x)\leq\frac{\sqrt3}{18}$ and equality holds only if $x=\frac{3-\sqrt3}{3}$.
\end{fact}
\emph{Proof of Fact \ref{aaa}.} Since $f'(x)=\frac{3x^2-6x+2}{4},$
then $f(x)$ is increasing when $x\in [0, \frac{3-\sqrt3}{3}]$ and decreasing when $x\in [\frac{3-\sqrt3}{3}, 1]$. Therefore $f(x)\leq f(\frac{3-\sqrt3}{3})=\frac{\sqrt3}{18}.$ \q
\begin{lemma}\label{B2n}
$\lambda(B(2, n-2))\leq\frac{\sqrt 3}{18}$ and $\underset{n \to +\infty}{\lim}\lambda(B(2, n-2))={\sqrt 3\over 18}.$
\end{lemma}
\emph{Proof of Lemma \ref{B2n}.} Let $\vec{x}=\{x_1,x_2,\dots,x_n\}$ be an optimum vector of $\lambda(B(2, n-2))$. Let $x_1+x_2=a$ and $b=1-a$. Then
\begin{eqnarray*}
\lambda(B(2, n-2))&\leq&\frac{a^2(1-a)}{4}+a\bigg(\frac{1-a}{n-2}\bigg)^2{n-2\choose 2}\\
&\leq&\frac{a^2(1-a)}{4}+\frac{a(1-a)^2}{2}.
\end{eqnarray*}
By Fact \ref{aaa}, $\lambda(B(2, n-2))\leq\frac{\sqrt 3}{18}$. Note that $'='$ holds only if $a=\frac{3-\sqrt 3}{3}$ and $n\to\infty$.
\q

\subsection{Preliminaries for the main Lemmas}
In this section we introduce two hypergraphs and show their Lagrangian are less than ${\sqrt3\over 18}$. When giving some proofs in Section \ref{sec4}, we will change our hypergraph by replacing some edges such that the new hypergraph has non-decreasing Lagrangian and is isomorphic with the subgraph of the following two hypergraphs.

$\mathbf{H_1:}$ the 3-graph with vertex set $[n]$ and edge set $E(H_1)=E(B(2, n-2))\setminus\{2ij: i,j\in [n]\setminus\{1, 2, 3, 4, 5, 6\}\}\cup\{345, 346\}$.

$\mathbf{H_2:}$ a 3-graph with vertex set $[n]$ and edge set $E(H_2)=E(B(2, n-2))\setminus\{2ij: i,j\in D\}\cup\{34i: i\in D\}$, where $D$ is a subset of $[n]\setminus \{1, 2, 3, 4\}$ with $|D|\geq 2$.

\begin{lemma}\label{h1}
$\lambda(H_1)<\frac{\sqrt3}{18}$.
\end{lemma}
\emph{Proof of Lemma \ref{h1}.} Let $\vec{x}=(x_1, x_2, ..., x_n)$ be an optimum vector of $\lambda(H_1)$. Assume that $x_1=a$, $x_2=b$, $x_3+x_4=c$, $x_5+x_6=d$, $e=1-a-b-c-d$. By Fact \ref{symmetry}, we may assume that $x_5=x_6=\frac{d}{2}$ and $x_3=x_4=\frac{c}{2}$. By Fact \ref{*}, we may assume that $a\geq b$. If $c=0$ or $d=0$, then $\lambda(H_1)<\lambda(B(2, n-2))\leq\frac{\sqrt 3}{18}$. If $b=0$, then
\begin{eqnarray*}
\lambda(H_1)&\leq&a(\frac{c^2}{4}+\frac{d^2}{4}+\frac{e^2}{2}+cd+ce+de)+\frac{c^2d}{4}\\
&\leq&\frac{a(c+d+e)^2}{2}+\frac{(\frac{c}{2}+\frac{c}{2}+d+e)^3}{27}=\frac{a(1-a)^2}{2}+\frac{(1-a)^3}{27}=f(a).\\
f'(a)&=&\frac{25a^2}{18}-\frac{16a}{9}+\frac{7}{18}.
\end{eqnarray*}
Note that $f(a)$ is increasing in $[0, \frac{7}{25}]$ and decreasing in $[{7\over 25}, 1]$, then $\lambda\leq f(\frac{7}{25})<\frac{\sqrt 3}{18}$. So we may assume that $a, b, c, d>0$.
By Fact \ref{fact2}, $$\frac{\partial\lambda}{\partial x_1}+\frac{\partial\lambda}{\partial x_2}=\frac{\partial\lambda}{\partial x_5}+\frac{\partial\lambda}{\partial x_6}.$$
Note that
$$\frac{\partial\lambda}{\partial x_1}+\frac{\partial\lambda}{\partial x_2}\geq bc+bd+be+cd+ce+de+\frac{c^2}{4}+\frac{d^2}{4}+ac+ad+ae+cd+ce+de+\frac{c^2}{4}+\frac{d^2}{4},$$
and $$\frac{\partial\lambda}{\partial x_5}+\frac{\partial\lambda}{\partial x_6}= ab+ac+ae+bc+be+(a+b)\frac{d}{2}+\frac{c^2}{4}+ab+ac+ae+bc+be+(a+b)\frac{d}{2}+\frac{c^2}{4}.$$
Therefore $$2cd+2ce+2de+\frac{d^2}{2}\leq 2ab+ac+ae+bc+be.$$
If $a+b\leq d$, then $2ab\leq\frac{d^2}{2}$, $ac+bc\leq cd$ and $ae+be\leq de$. So we have $a=b$ and $e=c=0$, a contradiction to $c>0$.
So we may assume that $a+b>d$, then
\begin{eqnarray*}
\lambda(H_1)&\leq&ab(c+d+e)+a(\frac{c^2}{4}+\frac{d^2}{4}+\frac{e^2}{2}+cd+ce+de)+b(\frac{c^2}{4}+\frac{d^2}{4}+cd+ce+de)+\frac{c^2d}{4}\\
&\leq&ab(1-a-b)+\frac{(a+b)(c+d+e)^2}{2}+\frac{c^2d}{4}-\frac{c^2a}{4}-\frac{c^2b}{4}\\
&<&ab(1-a-b)+\frac{(a+b)(1-a-b)^2}{2}\\
&\leq&\frac{(a+b)^2(1-a-b)}{4}+\frac{(a+b)(1-a-b)^2}{2}.
\end{eqnarray*}
By Fact \ref{aaa}, then $\lambda(H_1)<\frac{\sqrt 3}{18}.$
\q

\begin{lemma}\label{h2}
$\lambda(H_2)\leq\frac{\sqrt3}{18}$.
\end{lemma}
\emph{Proof of Lemma \ref{h2}.} Let $\vec{x}=(x_1, x_2, ..., x_n)$ be an optimum vector of $\lambda(H_2)$. Let $x_1=a$, $x_2=b$, $x_3+x_4=c$, $\sum_{v\in D}x_v=d$, $\sum_{v\in E}x_v=e$. We have
\begin{eqnarray*}
\lambda(H_2)%&\leq&ab(c+d+e)+a(\frac{c^2}{4}+\frac{d^2}{2}+\frac{e^2}{2}+cd+ce+de)+b(\frac{c^2}{4}+cd+ce+de)+\frac{c^2d}{4}\\
&\leq&ab(c+d+e)+a(\frac{c^2}{2}+\frac{d^2}{2}+\frac{e^2}{2}+cd+ce+de)+b(\frac{c^2}{4}+\frac{e^2}{2}+cd+ce+de)+\frac{c^2d}{4}\\
&=&\lambda(a, b, c, d, e)=\lambda
\end{eqnarray*}
under the constraint
\begin{equation}
\left\{
             \begin{array}{lr}
             a+b+c+d+e=1, &  \\
             a\geq 0, b\geq 0, c\geq 0, d\geq 0 \ and \ e\geq 0. &
             \end{array}
\right.
\end{equation}
Note that if $c=0$ or $d=0$, then
\begin{eqnarray*}
\lambda&<&\frac{(a+b)^2(1-a-b)}{4}+\frac{(a+b)(1-a-b)^2}{2},
\end{eqnarray*}
by Fact \ref{aaa}, then $\lambda\leq\frac{\sqrt 3}{18}$.
So we may assume that $c, d>0$. By Theorem \ref{KKT}, then $\frac{\partial \lambda}{\partial c}=\frac{\partial \lambda}{\partial d}$. By direct calculation,
\begin{eqnarray*}
\frac{\partial \lambda}{\partial c}&=&ab+a(c+d+e)+b(\frac{c}{2}+d+e)+\frac{cd}{2}\\
\frac{\partial \lambda}{\partial d}&=&ab+a(c+d+e)+b(c+e)+\frac{c^2}{4},
\end{eqnarray*}
then $\frac{c}{2}(\frac{c}{2}+b)=d(\frac{c}{2}+b)$, so $c=2d.$

By Fact \ref{*}, we may assume that $a\geq b$. We claim that $b>0$. If $b=0$, then
\begin{eqnarray*}
\lambda&=&a(\frac{c^2}{2}+\frac{d^2}{2}+\frac{e^2}{2}+cd+ce+de)+\frac{c^2d}{4}\\
&\leq&\frac{a(c+d+e)^2}{2}+\bigg(\frac{\frac{c}{2}+\frac{c}{2}+d+e}{3}\bigg)^3\\
&=&\frac{a(1-a)^2}{2}+\frac{(1-a)^3}{27}=f(a),\\
f'(a)&=&\frac{25a^2}{18}-\frac{16a}{9}+\frac{7}{18}.
\end{eqnarray*}
Note that $f(a)$ is increasing in $[0, \frac{7}{25}]$ and decreasing in $[\frac{7}{25}, 1]$, then $\lambda\leq f(\frac{7}{25})\leq\frac{\sqrt 3}{18}$. So we may assume that $a, b>0$. By Theorem \ref{KKT}, then $\frac{\partial \lambda}{\partial a}=\frac{\partial \lambda}{\partial b}$. Therefore $(a-b)(c+d+e)=\frac{c^2}{4}+\frac{d^2}{2}.$

We claim that $e>0$. If $e=0$, recall that we have shown that $c=2d$ and $(a-b)(c+d)=\frac{c^2}{4}+\frac{d^2}{2}$, then $a=b+\frac{d}{2}$. Since $a+b+c+d=1$, then $a=\frac{1}{2}-\frac{5d}{4}$ and $b=\frac{1}{2}-\frac{7d}{4}$. So
\begin{eqnarray*}
\lambda&=&\frac{-53d^3-12d^2+12d}{16}\\
\lambda'&=&\frac{-159d^2}{16}-\frac{3d}{2}+\frac{3}{4}.
\end{eqnarray*}
Note that $\lambda$ is increasing in $[0, \frac{2\sqrt{57}-4}{53}]$ and decreasing in $[\frac{2\sqrt{57}-4}{53}, 1]$, so $\lambda\leq 0.094.$

Therefore we may assume that $a, b, c, d, e>0$. Recall that $c=2d$. By Theorem \ref{KKT}, then $\frac{\partial \lambda}{\partial d}=\frac{\partial \lambda}{\partial e}$. Since
\begin{eqnarray*}
\frac{\partial \lambda}{\partial d}&=&ab+a(c+d+e)+b(c+e)+\frac{c^2}{4}\\
\frac{\partial \lambda}{\partial e}&=&ab+a(c+d+e)+b(c+d+e),
\end{eqnarray*}
then $d=b$. So
\begin{eqnarray*}
\lambda&<&\frac{(a+b)^2(1-a-b)}{4}+\frac{(a+b)(1-a-b)^2}{2}+\frac{c^2d}{4}-\frac{bc^2}{4}-\frac{bd^2}{2}\\
&<&\frac{(a+b)^2(1-a-b)}{4}+\frac{(a+b)(1-a-b)^2}{2}.
\end{eqnarray*}
By Fact \ref{aaa}, then $\lambda\leq\frac{\sqrt 3}{18}.$
\q

%%%%%%%%%%%%%%%%%%%%%%%%%%%%%%%%%%%%%%%%%%%%%%%%%%%%%

\section{Proofs of the main Lemmas}\label{sec4}
To complete the proof of Theorem \ref{main}, what remains is to show Lemma \ref{X2} to \ref{Y2}. In this section, we prove these lemmas. Throughout this section, let $G$ be a dense $K_4^3\cup e$-free 3-graph on vertex set $[n]$ with Lagrangian $\lambda(G)>\frac{\sqrt3}{18}$. Let $\vec x=\{x_1, x_2,\dots, x_n\}$ be an optimum vector for $\lambda(G)$. Since $\lambda(K_6^3)=\frac{5}{54}\leq\frac{\sqrt 3}{18}$, then  $n\geq 7$.

\begin{lemma}\label{TK43}\label{K43}(\cite{Baber},\cite{Raz})
$\pi(K_4^3)\leq 0.5615$.
\end{lemma}

By Fact \ref{dense}, Proposition \ref{relationlt} and Lemma \ref{K43}, $\pi_{\lambda}(K_4^3)=\pi(K_4^3)\leq 0.5615<3!\frac{\sqrt 3}{18}<3!\lambda(G)$. So $K_4^3\subseteq G$. Without loss of generality, assume that $\{1, 2, 3, 4\}$ forms a $K_4^3$ in $G$.

For $x, y \in V(G)$, let $N^*(x, y)=\{v: vxy\in E(G)\}$.

\begin{claim}\label{K53}
$K_5^3\nsubseteq G$.
\end{claim}
{\em Proof of Claim \ref{K53}.} Assume that $K_5^3\subseteq G$. Since $n\geq 7$, then there are $x, y\in [n]\setminus V(K_5^3)$, since $G$ is dense, then there exists $z\in N^*(x, y)$. The edge $xyz$ together with the $K_5^3$ contains a $K_4^3\cup e$, a contradiction.
\q

\begin{claim}\label{clique}
For $v\in [n]$, $\omega(G_v)\geq 3$ and $x_v<1-\sqrt {\frac{\sqrt 3}{3}\frac{\omega(G_v)}{\omega(G_v)-1}}$. Furthermore, if $v\in [n]\setminus\{1, 2, 3, 4\}$, then $3\leq\omega(G_v)\leq 4$. Therefore $x_v<0.0694$ if $\omega(G_v)=3$, and $x_v<0.12262$ if $\omega(G_v)=4$.
\end{claim}
{\em Proof of Claim \ref{clique}.} For any $v\in [n]$, applying Fact \ref{fact2} and Theorem \ref{MStheo}, we have
$$\frac{\sqrt 3}{6}<3\lambda=\frac{\partial\lambda}{\partial x_v}\leq\bigg(\frac{1-x_v}{\omega(G_v)}\bigg)^2{\omega(G_v)\choose 2}.$$
Then $$x_v<1-\sqrt {\frac{\sqrt 3}{3}\frac{\omega(G_v)}{\omega(G_v)-1}}.$$

If $v\in\{1, 2, 3, 4\}$, then it's clear that $\omega(G_v)\geq 3$. Let $v\notin [n]\setminus\{1, 2, 3, 4\}$. If $\omega(G_v)\leq 2$, then $x_v<0$, a contradiction. Therefore $\omega(G_v)\geq 3.$ If $\omega(G_v)\geq 5$, then $\omega(G_v)$ contains at least 1 vertex in $[n]\setminus\{1, 2, 3, 4\}$, therefore $G$ contains either a $K_5^3$ (if $\omega(G_v)$ contains only 1 vertex in $[n]\setminus\{1, 2, 3, 4\}$) or a $K_4^3\cup e$ (if $\omega(G_v)$ contains at least 2 vertices in $[n]\setminus\{1, 2, 3, 4\}$), a contradiction. Therefore $3\leq\omega(G_v)\leq 4$ for all $v\in [n]\setminus \{1, 2, 3, 4\}$.
\q

\begin{claim}\label{2}
For 2 vertices $a$ and $b$ in $[n]$, $x_a+x_b\leq\frac{3-\sqrt 3}{3}$.
\end{claim}
{\em Proof of Claim \ref{2}.} Applying Fact \ref{fact2}, we have
\begin{eqnarray*}
\frac{\sqrt 3}{3}&\leq&6\lambda(G)=\frac{\partial \lambda(G)}{\partial x_a}+\frac{\partial \lambda(G)}{\partial x_b}\\
&\leq&x_b(1-x_a-x_b)+\bigg(\frac{1-x_a-x_b}{n-2}\bigg)^2{{n-2}\choose 2}+x_a(1-x_a-x_b)+\bigg(\frac{1-x_a-x_b}{n-2}\bigg)^2{{n-2}\choose 2}\\
&\leq&(x_a+x_b)(1-x_a-x_b)+(1-x_a-x_b)^2\\
&=&1-(x_a+x_b).
\end{eqnarray*}
Therefore $x_a+x_b\leq\frac{3-\sqrt 3}{3}.$
\q

\begin{claim}\label{3}
If $G-\{v\}$ is $K_4^3$-free for some $v\in [n]$, then $x_v>0.0848$.
\end{claim}
{\em Proof of Claim \ref{3}.} Applying Lemma \ref{K43} and Proposition \ref{relationlt}, we have $$\lambda(G-\{v\}, \vec x)\leq(1-x_v)^3\frac{0.5615}{6}.$$
Therefore
\begin{eqnarray*}
\lambda(G)&\leq&(1-x_v)^3\frac{0.5615}{6}+x_v\frac{\partial \lambda(G)}{\partial x_v}\\
&=&(1-x_v)^3\frac{0.5615}{6}+3x_v\lambda(G),
\end{eqnarray*}
the last equality follows from Fact \ref{fact2}. So
\begin{eqnarray*}
\lambda(G)&\leq&\frac{0.5615}{6}\frac{(1-x_v)^3}{1-3x_v}.
\end{eqnarray*}
Note that $\frac{(1-x_v)^3}{1-3x_v}$ is increasing in $[0, \frac{1}{3})$. If $x_v\leq 0.0848$, then $\lambda(G)\leq0.09622\leq\frac{\sqrt 3}{18}$, a contradiction.
\q

Let $M^r_t$ be an $r$-graph on $tr$ vertices with $t$ disjoint edges.

\begin{claim}\label{8}
$n\geq 8$.
\end{claim}
{\em Proof of Claim \ref{8}.} Assume that and $V(G)=[7]$. Recall that $\{1, 2, 3, 4\}$ forms $K_4^3$. By Claim \ref{clique}, then $\omega(G_5), \omega(G_6), \omega(G_7)\leq 4$ and $x_5+x_6+x_7\leq 3\times 0.12262$. Therefore $x_1+x_2+x_3+x_4\geq 0.63214>4\times 0.158,$ then, without loss of generality, let $x_1\geq 0.158$. By Claim \ref{clique}, $x_v<0.151$ if $\omega(G_v)\leq 5$. Therefore $\omega(G_1)=6$. Since $G$ is $K_4^3\cup e$-free, then $G-\{1\}$ is $M_2^3$-free. Hefetz and Keevash (\cite{HK}) proved that $\pi_{\lambda}(M_2^3)\leq\frac{12}{25}$. So $\lambda(G-\{1\})\leq\frac{2}{25}$. Therefore
\begin{eqnarray*}
\lambda(G)&\leq&x_1\bigg(\frac{1-x_1}{6}\bigg)^2{6\choose 2}+\frac{2}{25}(1-x_1)^3\\
&=&\frac{5}{12}x_1(1-x_1)^2+\frac{2}{25}(1-x_1)^3=f(x_1).\\
f'(x_1)&=&\frac{5}{12}(1-x_1)^2-\frac{5}{6}x_1(1-x_1)-\frac{6}{25}(1-x_1)^2\\
&=&\frac{(1-x_1)(53-303x_1)}{300}.
\end{eqnarray*}
So $f(x_1)$ is increasing in $[0, \frac{53}{303}]$ and decreasing in $[\frac{53}{303}, 1]$, then $f(x_1)\leq f(\frac{53}{303})<0.095<\frac{\sqrt 3}{18}.$
\q

\subsection{$G$ is $K_5^{3-}$-free}
\begin{lemma}\label{K53-}
$G$ is $K_5^{3-}$-free.
\end{lemma}
{\em Proof of Lemma \ref{K53-}.} Assume that $K_5^{3-}\subseteq G$ with vertex set $\{1, 2, 3, 4, 5\}$ and $345\notin E(G)$. Since $G$ is $K_4^3\cup e$-free, then for any $x, y\in [n]\setminus\{1, 2, 3, 4, 5\}$, we have $N^*(x, y)\subseteq\{1, 2\}$. By Claim \ref{8}, $n\geq 8$. Recall that $\vec{x}=(x_1,x_2,\dots,x_n)$ is an optimal vector for $G$.

{\em Case 1.} $x_1, x_2\geq x_3, x_4, x_5$.

We say that a vertex $v\in [n]\setminus\{1, 2, 3, 4, 5\}$ is a {\em good vertex} if for the set of edges in $$B_v=\{vij\in E(G): \{i, j\}\in\{3, 4, 5\}^{(2)}\},$$ there exist the same number of triples in $$A_v=\{vij\in E(G^c): \{i, j\}\cap\{1, 2\}\not=\emptyset, \{i, j\}\subset\{1, 2, 3, 4, 5\}\}$$ such that $\sum_{vij\in B_v}x_vx_ix_j\leq\sum_{vij\in A_v}x_vx_ix_j$. In this case, we say that $B_v$ can be replaced by $A_v$. Otherwise, we call $v$ a {\em bad vertex}.

We call $v34, v35, v45$ bad edges for $v\in [n]\setminus\{1, 2, 3, 4, 5\}$. Note that for a good vertex $v\in [n]\setminus\{1, 2, 3, 4, 5\}$, replacing $B_v$ by $A_v$ in $G$ does not reduce the Lagrangian. Let $B$ be the set of all {\em bad vertices}. If $B=\emptyset$, then we can replace $B_v$ by $A_v$ in $G$ for each $v\in [n]\setminus\{1, 2, 3, 4, 5\}$, obtain $G^0$, and all edges in $E(G^0)$ are incident to 1 or 2. So $G^0\subseteq B(2, n-2)$. Hence $\lambda(G)\leq\lambda(G^0)\leq\lambda(B(2, n-2))=\frac{\sqrt 3}{18}$, a contradiction.

{\em Case 1.1.} There exists $v\in B$ such that exactly 1 of $v34, v35, v45$ is in $E(G)$.

Without loss of generality, let $v34\in E(G)$. Since $G$ is $K_5^3$-free, then at least 1 of $v12, v13, v23, v14, v24\notin E(G)$. Since $x_1, x_2\geq x_3, x_4$, then $v34$ can be replaced by that missing edge, so $v\notin B$, a contradiction.

{\em Case 1.2.} There exists $v\in B$ such that exactly 2 of $v34, v35, v45$ in $E(G)$.

Without loss of generality, let $v34, v35\in E(G)$. Since $\{v, 1, 2, 3, 4\}$ can't form a $K_5^3$, then at least 1 of $v12, v13, v23, v14, v24$ is not in $E(G)$.

If $v12\notin E(G)$, then $v13, v14, v23, v24\in E(G)$. Otherwise we can replace $v34, v35$ by that missing edge and $v12$ with $\lambda(G)$ non-decreasing. Contradict to $v\in B$. Therefore $\{v, 1, 3, 4\}, \{v, 2, 3, 4\}$ form $K_4^3$. Since $G$ is $K_4^3\cup e$-free, then $N^*(x, y)=\emptyset$ for $x, y\in [n]\setminus\{v, 1, 2, 3, 4, 5\}$, a contradiction to that $G$ is dense.

If $v14\notin E(G)$ (or $v24\notin E(G)$), then $\{v, 1, 3, 5\}, \{v, 2, 3, 5\}$ form $K_4^3$. Otherwise we can replace $v34, v35$ by $v14$ and one of the missing edges, a contradiction to $v\in B$. Since $G$ is $K_4^3\cup e$-free, then $N^*(x, y)=\emptyset$ for $x, y\in [n]\setminus\{v, 1, 2, 3, 4, 5\}$, a contradiction to that $G$ is dense.

If $v23\notin E(G)$ (or $v13\notin E(G)$), then $\{v, 1, 3, 4\}, \{v, 1, 3, 5\}, \{v, 1, 2, 4\}, \{v, 1, 2, 5\}$ form $K_4^3$. Since $G$ is $K_4^3\cup e$-free, then for all $x, y\in [n]\setminus\{v, 1, 2, 3, 4, 5\}$, we have $N^*(x, y)=\{1\}$ and $x34, x35\notin E(G)$. So for any $x\in [n]\setminus\{v, 1, 2, 3, 4, 5\}$ only, possibly, $x45\in E(G)$ is a bad edge incident to $x$. By the proof of {\em Case 1.1}, $x\notin B$, so $B=\{v\}$. Replacing $v35$ by $v23,$ adding $345$ and replacing $B_x$ by $A_x$ for all $x\in [n]\setminus\{v, 1, 2, 3, 4, 5\}$, we obtain $G^0$. Note that $G^0$ is contained in an isomorphic copy of $H_1$ (view $v$ in $G^0$ as 6 in $H_1$), then $\lambda(G)\leq\lambda(G^0)\leq\lambda(H_1)\leq\frac{\sqrt 3}{18}$ by Lemma \ref{h1}, a contradiction to $\lambda(G)>\frac{\sqrt 3}{18}.$

{\em Case 1.3.} There exists $v\in B$ such that $v34, v35, v45\in E(G)$.

Since $v\in B$, then at most two of $\{v12, v13, v23, v14, v24, v15, v25\}$ are not in $E(G)$, otherwise we can replace $B_v$ by the three missing edges. We claim that there are 2 of those edges not in $E(G)$. Suppose that there is only 1 of those edges missing in $E(G)$. If only $v12\notin E(G)$, then $\{v, 1, 3, 4\}, \{v, 2, 3, 4\}$ form $K_4^3$. Since $G$ is $K_4^3\cup e$-free and $|G|\geq 8$, then there are $x, y\in [n]\setminus\{v, 1, 2, 3, 4, 5\}$ such that $N^*(x, y)=\emptyset$, a contradiction to that $G$ is dense. If  only $vij\notin E(G)$ for some $i\in\{1, 2\}$ and some $j\in\{3, 4, 5\}$,
then $\{v, 1, 2, k, l\}$ form a $K_5^3$, where $\{k, l\}=\{3, 4, 5\}\setminus\{j\}$, a contradiction to Claim \ref{K53}. So we may assume that there are two of those edges not in $E(G)$.

If $v12\notin E(G)$ and $vij\notin E(G)$, where $i\in\{1, 2\}$ and $j\in\{3, 4, 5\}$, then $\{v, 1, k, l\}$ and $\{v, 2, k, l\}$ ($\{k, l\}=\{3, 4, 5\}\setminus\{j\}$) form $K_4^3$. Since $G$ is $K_4^3\cup e$-free, then $N^*(x, y)=\emptyset$ for $x, y\in [n]\setminus\{v, 1, 2, 3, 4, 5\}$, contradicting to $G$ being dense.

If $v1i, v2j\notin E(G)$, where $i, j\in \{3, 4, 5\}$, then $\{v, 1, p, q\}$ and $\{v, 2, s, t\}$ form $K_4^3$, where $\{p, q\}=\{3, 4, 5\}\setminus\{i\}$ and $\{s, t\}=\{3, 4, 5\}\setminus\{j\}$. Since $G$ is $K_4^3\cup e$-free, then $N^*(x, y)=\emptyset$ for $x, y\in [n]\setminus\{v, 1, 2, 3, 4, 5\}$, contradicting to $G$ being dense.
%f $v1i, v2i\notin E(G)$, where $i\in \{3, 4, 5\}$), then $\{v, 1, 4, 5\}$ and $\{v, 2, 4, 5\}$ form $K_4^3$. Since $G$ is $K_4^3\cup\{e\}$-free, then $N^*(x, y)=\emptyset$ for $x, y\in [n]\setminus\{v, 1, 2, 3, 4, 5\}$, contradicting to $G$ being dense.

If $vij, vik\notin E(G)$, where $i\in\{1, 2\}$ and $\{j, k\}\in \{3, 4, 5\}^{(2)}$, without loss of generality, assume that $v13, v14\notin E(G)$, then $\{v, 2, 3, 4\}, \{v, 2, 3, 5\}, \{v, 2, 4, 5\}, \{v, 1, 2, 5\}$ form $K_4^3$. Since $G$ is $K_4^3\cup e$-free, then for all $x, y\in [n]\setminus\{v, 1, 2, 3, 4, 5\}$, we have $N^*(x, y)=\{2\}$ and $x34\notin E(G)$. By the proof of {\em Case 1.2}, $x\notin B$, so $B=\{v\}$. Replacing $B_x$ by $A_x$ for all $x\in [n]\setminus\{v, 1, 2, 3, 4, 5\}$, replacing $v35, v45$ by $v13, v14$ and adding 345, we obtain $G^0$ which is contained in an isomorphic copy of $H_1$ (view $v$ in $G^0$ as 6 in $H_1$). So $\lambda(G)\leq\lambda(G^0)\leq\lambda(H_1)\leq\frac{\sqrt 3}{18}$ by Lemma \ref{h1}, a contradiction to $\lambda(G)>\frac{\sqrt 3}{18}.$

\bigskip

{\em Case 2.} $x_1, x_3\geq x_2, x_4, x_5.$

A vertex $v\in [n]\setminus\{1, 2, 3, 4, 5\}$ is a {\em good vertex} if for the edges in
$$B_v=\{vij\in E(G): ij\in\{2, 4, 5\}^{(2)}\},$$ there exist the same number of triples in
$$A_v=\{vi'j'\in E(G^c): \{i', j'\}\cap\{1, 3\}\not=\emptyset, \{i', j'\}\subset\{1, 2, 3, 4, 5\}\}$$ such that the substitute $vi'j'$ for $vij\in B_v$ satisfies $|\{i', j'\}\cap\{i, j\}|=1$ or $\{i', j'\}=\{1, 3\}$. Note that $\sum_{vij\in B_v}x_vx_ix_j\leq\sum_{vij\in A_v}x_vx_ix_j$. In this case, we say that $B_v$ can be replaced by $A_v$. Otherwise, we call $v$ a {\em bad vertex}. We call $v24, v25, v45$ {\em bad edges} for $v\in [n]\setminus\{1, 2, 3, 4, 5\}$.

Let $B$ be the vertex set containing all {\em bad vertices}.

{\em Case 2.1.} There exists $v\in B$ and exactly one of $\{v24, v25, v45\}$ is in $E(G)$.

{\em Case 2.1.1.} $v45\in E(G)$.

Since $v\in B$, then $v13, v14, v34, v15, v35\in E(G)$. Therefore $\{v, 1, 3, 4\}, \{v, 1, 3, 5\}, \{v, 1, 4, 5\}$ form $K_4^3$, then for all $x, y\in [n]\setminus\{v, 1, 2, 3, 4, 5\}$ we have $N^*(x, y)=\{1\}$ and $x23, x24, x25\notin E(G)$. Otherwise $K_4^3\cup e\subseteq G$. For any $u\in B$, since $u24, u25\notin E(G)$, then $u45\in E(G)$ and $\{u, 1, 3, 4\}, \{u, 1, 3, 5\},\\ \{u, 1, 4, 5\}$ form $K_4^3$.

If $B=\{v\}$, then let $E=[n]\setminus\{1, 2, 3, 4, 5, v\}$. Note that $N^*(x, y)=1$ for $x, y\in E$. Replacing $B_x$ by $A_x$ for $x\in E$, we obtain $G^0$. Then  $G^0$ is contained in an isomorphic copy of $H_1$ (view 3 in $G^0$ as 2 in $H_1$, view $245, v45$ in $G^0$ as $345, 346$ in $H_1$). So $\lambda(G)\leq\lambda(G^0)\leq\lambda(H_1)\leq\frac{\sqrt 3}{18}$ by Lemma \ref{h1}, a contradiction.

If $B=\{v, u\}$, then let $E=[n]\setminus(B\cup\{1, 3, 4, 5\})$. Note that $2\in E$. Since $x$ is a good vertex for $x\in E\setminus\{2\}$, then $x45$ can be replaced by 1 of $\{x13, x14, x15, x34, x35\}$. Note that $N^*(x, y)=1$ for $x, y\in E\cup\{2\}$. Replace $uv2$ by $uv3$, then the obtained $G^0$ is contained in an isomorphic copy of $H_1$ (view 3 in $G^0$ as 2 in $H_1$, view $v45, u45$ in $G^0$ as $345, 346$ in $H_1$). So $\lambda(G)\leq\lambda(G^0)\leq\lambda(H_1)\leq\frac{\sqrt 3}{18}$ by Lemma \ref{h1}, a contradiction.

If $|B|\geq 3$, then $N^*(x, y)=\{1\}$ for $x, y\in B$. Otherwise if $xy2\in E(G)$, then $\{z, 1, 4, 5\}\cup\{x, y, 2\}$ forms $K_4^3\cup e$ in $G$ for $x, y, z\in B$. Let $E=[n]\setminus(B\cup\{1, 3, 4, 5\})$. Replacing $B_x$ by $A_x$ for each $x\in E$, we obtain $G^0$ which is contained in an isomorphic copy of $H_2$ (view $B$ in $G^0$ as $D$ in $H_2$, view 3 in $G^0$ as 2 in $H_2$, view $i45 (i\in B)$ in $G^0$ as $i34 (i\in D)$ in $H_2$). So $\lambda(G)\leq\lambda(G^0)\leq\lambda(H_2)\leq\frac{\sqrt 3}{18}$ by Lemma \ref{h2}, a contradiction.

{\em Case 2.1.2.} $v2i\in E(G)$, where $i$ is 4 or 5.

Since $v\in B$, then $v12, v23, v1i, v3i, v13\in E(G)$, so $\{v, 1, 2, 3, i\}$ forms a $K_5^3$, a contradiction to Claim \ref{K53}.

{\em Case 2.2.} There exists $v\in B$ and exactly two of $\{v24, v25, v45\}$ are in $E(G)$.

{\em Case 2.2.1.} $v24, v25\in E(G)$.

Since $\{v, 1, 2, 3, 4\}$ can't form a $K_5^3$, then at least one of $\{v12, v23, v14, v34, v13\}$ is not in $E(G)$.

If $v12\notin E(G)$, then $v13, v14, v23, v34\in E(G)$. Otherwise we can replace $v24, v25$ by $v12$ and that missing edge, a contradiction to $v\notin B$. Therefore $\{v, 1, 3, 4\}, \{v, 2, 3, 4\}$ form $K_4^3$. Since $G$ is $K_4^3\cup e$-free, then $N^*(x, y)=\emptyset$ for all $x, y\in [n]\setminus\{v, 1, 2, 3, 4, 5\}$, a contradiction to that $G$ is dense.

If $v23\notin E(G)$, then $v12, v13, v14, v15, v34, v35\in E(G)$. Therefore $\{v, 1, 2, 4\}, \{v, 1, 2, 5\}, \{v, 1, 3, 4\},\\\{v, 1, 3, 5\}$ form $K_4^3$, then for all $x, y\in [n]\setminus\{v, 1, 2, 3, 4, 5\}$ we have $N^*(x, y)=\{1\}$ and $x24, x25\notin E(G)$,  by the proof of {\em Case 2.1}, $x\notin B$, so $B=\{v\}$. Let $E=[n]\setminus\{1, 2, 3, 4, 5, v\}$. Deleting $v25$, adding $v23$, and replacing $B_x$ by $A_x$ for $x\in E$, we obtain $G^0$ which is contained in an isomorphic copy of $H_1$ (view 3, 2 in $G^0$ as 2, 3 in $H_1$ respectively, view $v24$ in $G^0$ as $346$ in $H_1$). So $\lambda(G)\leq\lambda(G^0)\leq\lambda(H_1)\leq\frac{\sqrt 3}{18}$ by Lemma \ref{h1}, a contradiction.

If $v34\notin E(G)$ (or $v14\notin E(G)$ or both), then $\{v12, v13, v15, v23, v35\}\subset E(G)$, so $\{v, 1, 3, 5\}, \{v, 2, 3, \\5\}$ form $K_4^3$. Since $G$ is $K_4^3\cup e$-free, then $N^*(x, y)=\emptyset$ for all $x, y\in [n]\setminus\{v, 1, 2, 3, 4, 5\}$, a contradiction to that $G$ is dense.

If $v13\notin E(G)$, then $\{v12, v14, v15, v23, v34, v35\}\subseteq E(G)$, otherwise we can replace $v24, v25$ by that missing edge and $v13$, so $\{v, 1, 2, 4\}, \{v, 1, 2, 5\}, \{v, 2, 3, 4\}, \{v, 2, 3, 5\}$ form $K_4^3$. Since $G$ is $K_4^3\cup e$-free, then $N^*(x, y)=\{2\}$ and $x34, x35, x14\notin E(G)$ for all $x, y\in [n]\setminus\{v, 1, 2, 3, 4, 5\}$. Note that $x\notin B$ for $x\in [n]\setminus\{v, 1, 2, 3, 4, 5\}$, otherwise we can replace the all bad edges incident to $x$ by $x34, x35, x14$. Deleting $v25$, adding $v13$, replacing $B_x$ by $A_x$ and replacing $xv2, xy2$ by $xv3, xy1$ for $x, y\in [n]\setminus\{1, 2, 3, 4, 5, v\}$, we obtain $G^0$ which is contained in an isomorphic copy of $H_1$ (view 3, 2 in $G^0$ as 2, 3 in $H_1$ respectively, view $v24$ in $G^0$ as $346$ in $H_1$). So $\lambda(G)\leq\lambda(G^0)\leq\lambda(H_1)\leq\frac{\sqrt 3}{18}$ by Lemma \ref{h1}, a contradiction.

{\em Case 2.2.2.} $v24, v45\in E(G)$ (the proof for $v25, v45\in E(G)$ is identical).

Since $\{v, 1, 2, 3, 4\}$ can't form a $K_5^3$, then at least 1 of $\{v23, v12, v13, v14, v34\}$ is not in $E(G)$.

If $v23\notin E(G)$, since $v\in B$, then $v13, v14, v15, v34, v35\in E(G)$. So $\{v, 1, 3, 4\}, \{v, 1, 3, 5\}, \{v, 1, 4, 5\}$ form $K_4^3.$ Since $G$ is $K_4^3\cup e$-free, then $N^*(x, y)=\{1\}$ and $x23, x24, x25\notin E(G)$ for all $x, y\in [n]\setminus\{v, 1, 2, 3, 4, 5\}$. By the proof of  {case 2.1}, then $x\notin B$, so $B=\{v\}$. Let $E=[n]\setminus\{1, 2, 3, 4, 5, v\}$. Deleting $v24$, adding $v23$, deleting $xv2$, adding $xv3$ for $x\in E$, and replacing $B_x$ by $A_x$ for $x\in E$, we obtain $G^0$ which is contained in an isomorphic copy of $H_1$ (view 3 in $G^0$ as 2 in $H_1$, view $\{245, v45\}$ in $G^0$ as $\{345, 346\}$ in $H_1$). So $\lambda(G)\leq\lambda(G^0)\leq\lambda(H_1)\leq\frac{\sqrt 3}{18}$ by Lemma \ref{h1}, a contradiction. So we may assume that $v23\in E(G)$.

If $v12\notin E(G)$, since $v\in B$, then $v13, v14, v23, v34\in E(G)$, so $\{v, 1, 3, 4\}, \{v, 2, 3, 4\}$ form $K_4^3$. Since $G$ is $K_4^3\cup e$-free, then $N^*(x, y)=\emptyset$ for $x, y\in [n]\setminus\{v, 1, 2, 3, 4, 5\}$, contradicting to $G$ being dense. So we may assume that $v12\in E(G)$.

If $v13\notin E(G)$, since $v\in B$, then $v14, v15, v23, v34\in E(G)$, so $\{v, 1, 4, 5\}, \{v, 2, 3, 4\}$ form $K_4^3$. Since $G$ is $K_4^3\cup e$-free, then $N^*(x, y)=\emptyset$ for $x, y\in [n]\setminus\{v, 1, 2, 3, 4, 5\}$, contradicting to $G$ being dense. So we may assume that $v13\in E(G)$.

If $v14\notin E(G)$, since $v\in B$, then $v13, v15, v23, v34, v35\in E(G)$, so $\{v, 1, 3, 5\}, \{v, 2, 3, 4\}$ form $K_4^3$. Since $G$ is $K_4^3\cup e$-free, then $N^*(x, y)=\emptyset$ for $x, y\in [n]\setminus\{v, 1, 2, 3, 4, 5\}$, contradicting to $G$ being dense. So we may assume that $v14\in E(G)$.

If $v34\notin E(G)$, since $v\in B$, then $v12, v13, v14, v15, v23, v35\in E(G)$, so $\{v, 1, 2, 3\}, \{v, 1, 2, 4\}, \{v, 1, 3,\\ 5\}, \{v, 1, 4, 5\}$ form $K_4^3$. Since $G$ is $K_4^3\cup e$-free, then $N^*(x, y)=\{1\}$ and $x24, x45\notin E(G)$ for all $x, y\in [n]\setminus\{v, 1, 2, 3, 4, 5\}$, which means, by the proof of {\em Case 2.1.}, $x\notin B$ and $B=\{v\}$. Let $E=[n]\setminus\{1, 2, 3, 4, 5, v\}$. Deleting $v24$, adding $v34$, replacing $B_x$ by $A_x$, deleting $xv2$ and adding $xv3$ for $x\in E$, we obtain $G^0$ which is contained in an isomorphic copy of $H_1$ (view 3 in $G^0$ as 2 in $H_1$, view $\{245, v45\}$ in $G^0$ as $\{345, 346\}$ in $H_1$). So $\lambda(G)\leq\lambda(G^0)\leq\lambda(H_1)\leq\frac{\sqrt 3}{18}$ by Lemma \ref{h1}, a contradiction.

{\em Case 2.3.} There exists $v\in B$ and $v24, v25, v45\in E(G)$.

Since $v\in B$, then at most 2 of $\{v12, v13, v14, v15, v23, v34, v35\}$ are not in $E(G)$, otherwise we can replace $B_v=\{v24, v25, v45\}$ by those 3 missing edges in $A_v$.
We claim that exactly 2 of $\{v12, v13, v14, v15, v23, \\v34, v35\}$ are not in $E(G)$. Otherwise there is only 1 of those edges not in $E(G)$.
If only $v12\notin E(G)$ (or only $v23\notin E(G)$), then $\{v, 1, 4, 5\}$ , $\{v, 2, 4, 5\}$ form $K_4^3$. Since $G$ is $K_4^3\cup e$-free, then $N^*(x, y)=\emptyset$ for all $x, y\in [n]\setminus\{v, 1, 2, 3, 4, 5\}$, a contradiction to $G$ being dense.
If only $v13$ (or one of $\{v14, v15, v34, v35\})$ is not in $E(G)$, then $\{v, 1, 2, 4, 5\}$ forms an $K_5^3$, a contradiction. So we can assume that there are exactly two of $\{v12, v13, v14, v15, v23, v34, v35\}$  not in $E(G)$.

If $v12\notin E(G)$, and $v23\notin E(G)$ (the discussion for $v23$ replaced by 1 of $\{v13, v14, v15, v34, v35\}$ is similar), then $\{v, 1, 4, 5\}, \{v, 2, 4, 5\}$ form $K_4^3$. Since $G$ is $K_4^3\cup e$-free, then $N^*(x, y)=\emptyset$ for $x, y\in [n]\setminus\{v, 1, 2, 3, 4, 5\}$, contradicting to $G$ being dense. So we may assume that $v12\in E(G)$.

If $v13\notin E(G)$, and $v23\notin E(G)$ (or 1 of $\{v34, v35\}\notin E(G)$), then $\{v, 1, 4, 5\}, \{v, 2, 4, 5\}$ form $K_4^3$. Therefore $N^*(x, y)=\emptyset$ for all $x, y\in [n]\setminus\{v, 1, 2, 3, 4, 5\}$, a contradiction. If $v13\notin E(G)$ and $v14\notin E(G)$ (or $v15\notin E(G)$, since 4 and 5 are symmetric, we only discuss $v14$ here), then $\{v, 2, 3, 4\}, \{v, 2, 3, 5\}, \{v, 2, 4, 5\}, \{v, 1, 2, 5\}$ form $K_4^3$. Therefore $N^*(x, y)=\{2\}$ and $x13, x14, x15, xv3\notin E(G)$ for all $x, y\in [n]\setminus\{v, 1, 2, 3, 4, 5\}$.
Since we can replace $B_x$ by $x13, x14, x15$ for $x\in [n]\setminus\{v, 1, 2, 3, 4, 5\}$, then $x\notin B$, so $B=\{v\}$. Let $E=[n]\setminus\{1, 2, 3, 4, 5, v\}$. Deleting $v24$ and $v25$, adding $v13$ and $v14$, deleting $xv2$, adding $xv3$, deleting $xy2$, adding $xy1$ for $x, y\in E$, and replacing $B_x$ by $A_x$ for $x\in E$, we obtain $G^0$. Note that $G^0$ is contained in an isomorphic copy of $H_1$ (view 3 in $G^0$ as 2 in $H_1$, view $\{245, v45\}$ in $G^0$ as $\{345, 346\}$ in $H_1$). So $\lambda(G)\leq\lambda(G^0)\leq\lambda(H_1)\leq\frac{\sqrt 3}{18}$ by Lemma \ref{h1}, a contradiction. So we may assume that $v13\in E(G)$.
%If $v34\notin E(G)$ and $v14\notin E(G)$, then $\{v, 1, 3, 5\}, \{v, 2, 3, 5\}$ form $K_4^3$. Therefore $N^*(x, y)=\emptyset$ for all $x, y\in [n]\setminus\{v, 1, 2, 3, 4, 5\}$, a contradiction. If $v34\notin E(G)$ and $v23\notin E(G)$, then $\{v, 1, 4, 5\}, \{v, 2, 4, 5\}$ form $K_4^3$. Therefore $N^*(x, y)=\emptyset$ for all $x, y\in V(G)\setminus\{v, 1, 2, 3, 4, 5\}$, a contradiction. If $v34, v15\notin E(G)$, then $\{v, 1, 2, 3\}$, $\{v, 1, 2, 4\}$, $\{v, 2, 3, 5\}$ and $\{v, 2, 4, 5\}$ form $K_4^3$.
%Therefore $N^*(x, y)=\{2\}$ and $x13, x14, x35, x45\notin E(G)$ for all $x, y\in V(G)\setminus\{v, 1, 2, 3, 4, 5\}$. Note that $B=\{v\}$, otherwise we can replace $B_x$ by $x13, x14$ for $x\in V(G)\setminus\{v, 1, 2, 3, 4, 5\}$. Let $D=\{v, 2\}$ and $E=V(G)\setminus(D\cup\{1, 3, 4, 5\})$. Delete $v24$ and $v25$ add $v15$ and $v34$. Delete $xv2$ add $xv1$ and delete $xy2$ add $xy1$ for $x, y\in E$. Replace $B_x$ by $A_x$ for $x\in E$. We obtain $G^0$ and $G^0\subseteq H_1$. By Lemma \ref{h1}, then $\lambda(G)\leq\lambda(G^0)\leq\frac{\sqrt 3}{18}$, a contradiction. So we may assume that $v34\in E(G).$

If $v14\notin E(G)$, and $v23\notin E(G)$ (or $v34\notin E(G)$), then $\{v, 1, 3, 5\}$ and $\{v, 2, 4, 5\}$ form $K_4^3$. Therefore $N^*(x, y)=\emptyset$ for all $x, y\in [n]\setminus\{v, 1, 2, 3, 4, 5\}$, a contradiction. If $v14, v15\notin E(G)$ ($v14, v35\notin E(G)$ is similar), then $\{v, 1, 2, 3\}, \{v, 2, 3, 4\}, \{v, 2, 3, 5\}, \{v, 2, 4, 5\}$ form $K_4^3$. Therefore $N^*(x, y)=\{2\}$ and $x13, x14, x15, x45\notin E(G)$ for all $x, y\in [n]\setminus\{v, 1, 2, 3, 4, 5\}$. So we can replace $B_x$ by $x13, x14, x15$ for $x\in [n]\setminus\{v, 1, 2, 3, 4, 5\}$, then $x\notin B$,  and $B=\{v\}$. Let $E=[n]\setminus\{1, 2, 3, 4, 5, v\}$. Deleting $v24$ and $v25$, adding $v14$ and $v15$, deleting $xv2$, adding $xv3$, deleting $xy2$, adding $xy1$ for $x, y\in E$, and replacing $B_x$ by $A_x$ for $x\in E$, we obtain $G^0$ which is contained in an isomorphic copy of $H_1$ (view 3 in $G^0$ as 2 in $H_1$, view $\{245, v45\}$ in $G^0$ as $\{345, 346\}$ in $H_1$). By Lemma \ref{h1}, then $\lambda(G)\leq\lambda(G^0)\leq\frac{\sqrt 3}{18}$, a contradiction. So we may assume that $v14\in E(G)$.

Similar to $v14$, we may assume that $v15\in E(G)$, but then $\{v, 1, 2, 4, 5\}$ forms a $K_5^3$, a contradiction to Claim \ref{K53}.

\bigskip

{\em Case 3.} $x_3, x_4\geq x_1, x_2, x_5$.

A vertex $v\in [n]\setminus\{1, 2, 3, 4, 5\}$ is a {\em good vertex} if for the edges in $$B_v=\{vij\in E(G): ij\in\{1, 2, 5\}^{(2)}\},$$ there exist the same number of triples in $$A_v=\{vi'j'\in E(G^c): \{i', j'\}\cap\{3, 4\}\not=\emptyset, \{i', j'\}\subset\{1, 2, 3, 4, 5\}\}$$ such that the substitute $vi'j'$ for $vij\in B_v$ satisfies $|\{i', j'\}\cap\{i, j\}|=1$ or $\{i', j'\}=\{3, 4\}$. Note that $\sum_{vij\in B_v}x_vx_ix_j\leq\sum_{vij\in A_v}x_vx_ix_j$. In this case, we say that $B_v$ can be replaced by $A_v$. Otherwise we call $v$ a {\em bad vertex}. We call $v12, v15, v25$ {\em bad edges} for $v\in [n]\setminus\{1, 2, 3, 4, 5\}$.

Let $B$ be the set of all {\em bad vertices}. Let $E=[n]\setminus(B\cup\{1, 2, 3, 4, 5\})$.

{\em Case 3.1.} There exists $v\in B$ and there is exactly 1 of $v12, v15, v25\in E(G)$.

{\em Case 3.1.1.} $v12\in E(G)$.

Since $v\in B$, then $\{v13, v14, v23, v24, v34\}\subseteq E(G)$, so $\{v, 1, 2, 3, 4\}$ forms a $K_5^3$, a contradiction to Claim \ref{K53}.

{\em Case 3.1.2.} $v15\in E(G)$ (the case $v25\in E(G)$ is similar).

Since $v\in B$, then $v13, v14, v34, v35, v45\in E(G)$. So $\{v, 1, 3, 4\}, \{v, 1, 3, 5\}, \{v, 1, 4, 5\}$ form $K_4^3$. Since $G$ is $K_4^3\cup e$-free, then $N^*(x, y)=\{1\}$ and $v12, v25, x23, x24, x25\notin E(G)$ for all $x, y\in [n]\setminus\{v, 1, 2, 3, 4, 5\}$. Therefore for all $u\in B$, only possibly, $u12, u15\in E(G)$. Since we can replace $u12$ by $u23$ and $u\in B$, then $u15$ can't be replaced. So $\{u, 1, 3, 4\}, \{u, 1, 3, 5\}, \{u, 1, 4, 5\}$ form $K_4^3$.

If $|B|=1$, i.e. $B=\{v\}$, then replacing $B_x$ by $A_x$, deleting $xv1, xv2, xy1$, adding $xv3, xv4, xy3$ for all $x, y\in E$, we obtain $G^0$ which is contained in an isomorphic copy of $H_1$ (view 3, 4 in $G^0$ as 1, 2 in $H_1$ respectively, view $\{v15, 125\}$ in $G^0$ as $\{345, 346\}$ in $H_1$). Then $\lambda(G)\leq\lambda(G^0)\leq\frac{\sqrt 3}{18}$ by Lemma \ref{h1}, a contradiction.

If $|B|=2$, i.e. $B=\{v, v'\}$, then $v'15\in E(G)$ (or $v'12, v'15\in E(G)$ but we can replace $v'12$ by $v'23$). Deleting $vv'1, vv'2$, adding $vv'3, vv'4$, and replacing $xy1, xv1, xv'1, x12, 125$ by $xy3, xv3, xv'3, x23, 345$ respectively for all $x, y\in E$, we obtain $G^0$. View $3, 4, 1, 5, v, v'$ in $G^0$ as $1, 2, 3, 4, 5, 6$ in $H_1$, respectively. Note that $N^*(x, 2)=\{3\}$ in $G^0$ for $x\in E$, so $G^0$ is contained in an isomorphic copy of $H_1$. Hence $\lambda(G)\leq\lambda(G^0)\leq\frac{\sqrt 3}{18}$ by Lemma \ref{h1}, a contradiction.

If $|B|\geq3$, then $vv'2, xv2\notin E(G)$ for any $v, v'\in B$ and $x\in E$, since otherwise $\{v'', 1, 3, 4\}\cup\{v, v', 2\}$ forms $K_4^3\cup e$ for $v''\in B$ and $x\in E$. Replacing $vv'1$ by $vv'3$, deleting $xy1, xv1$, adding $xy3, xv3$ for all $x, y\in E$ and $v, v'\in B$, and replacing $B_x$ by $A_x$ for $x\in E$, we obtain $G^0$. So $N^*(v, v')=\{3\}$ for $v, v'\in B$ in $G^0$. View $\{3, 4\}$ in $G^0$ as $\{1, 2\}$ in $H_2$, view $v15$ in $G^0$ ($v\in B$) as $i34$ ($i\in D$) in $H_2$, then $G^0$ is contained in an isomorphic copy of $H_2$. So $\lambda(G)\leq\lambda(G^0)\leq\lambda(H_2)\leq\frac{\sqrt 3}{18}$ by Lemma \ref{h2}, a contradiction.

{\em Case 3.2.} $v\in B$ and exactly 2 of $v12, v15, v25\in E(G)$.

{\em Case 3.2.1.} $v12, v15\in E(G)$. (The case that $v12, v25\in E(G)$ is similar.)

Since $\{v, 1, 2, 3, 4\}$ can't form a $K_5^3$, then there is at least 1 of $\{v13, v14, v23, v24, v34\}$ not in $E(G)$. Since $v\in B$, then we may assume that $v35, v45\in E(G)$. Otherwise we may replace $v12, v15$ by 1 of $\{v13, v14, v23, v24, v34\}$ and 1 of $v35, v45\in E(G)$, a contradiction. If $v13\notin E(G)$ (or $v14\notin E(G)$), since $v\in B$, then $\{v, 2, 3, 4\}$ and $\{v, 1, 4, 5\}$ (or $\{v, 1, 3, 5\}$) form $K_4^3$, then $N^*(x, y)=\emptyset$ for all $x, y\in [n]\setminus\{v, 1, 2, 3, 4, 5\}$. So we may assume that $v13, v14\in E(G)$.

If $v23\notin E(G)$ (or $v24\notin E(G)$ or both), since $v\in B$, then $v34\in E(G)$. So $\{v, 1, 3, 4\}, \{v, 1, 3, 5\}, \{v, 1,\\ 4, 5\}$ form $K_4^3$. So $N^*(x, y)=\{1\}$ and $x23, x24, x25\notin E(G)$ for all $x, y\in [n]\setminus\{v, 1, 2, 3, 4, 5\}$. If $v'\in B\setminus\{v\}$, note that $v'25\notin E(G)$, only possibly, $v'12, v'15\in E(G)$. And $v'12$ can be replaced by $v'23$, so $v'15$ can't be replaced, therefore $\{v', 1, 3, 4\}, \{v', 1, 3, 5\}, \{v', 1, 4, 5\}$ form $K_4^3$.

If $B=\{v\}$, then replacing $B_x$ by $A_x$ for $x\in E$, deleting $xy1, xv1, xv2$, and adding $xy3, xv3, xv4$ for $x, y\in E$, we obtain $G^0$ which is contained in an isomorphic copy of $H_1$ (view $3, 4, v$ in $G^0$ as $1, 2, 6$ in $H_1$ respectively, view $\{v15, 125\}$ in $G^0$ as $\{345, 346\}$ in $H_1$), a contradiction.

If $B=\{v, v'\}$, then replace $v12$ by $v23$. Delete $125, 1vv', 2vv', xy1$ and add $345, 3vv', 4vv', xy3$ for $x, y\in E$, respectively. Replace $x12$ by $x23$. Since $x$ is a good vertex for $x\in E$, then $x15$ can be replaced by 1 of the missing edges in $\{x13, x14, x34, x35, x45\}$. Let $G^0$ be the resulting 3-graph, note that $N^*(x, y)=\{3\}$ for $x, y\in E\cup\{2\}$ in $G^0$ (view $3, 4, 1, 5, v, v'$ in $G^0$ as $1, 2, 3, 4, 5, 6$ in $H_1$ respectively). Then $G^0$ is contained in an isomorphic copy of $H_1$, a contradiction.

If $|B|\geq 3$, then replace $v12$ by $v23$ for $v\in B$. Delete $125$ and add $345$ for $x\in E$. Replace $B_x$ by $A_x$ for $x\in E$. Since $G$ is $K_4^3\cup e$-free, then $N^*(v, v')=\{1\}$ for any $v, v'\in B$. Replace $vv'1$ by $vv'3$ for $v, v'\in B$. Let $G^0$ be the resulting 3-graph. Note that $N^*(x, y)=\{3\}$ for $x, y\in B$ in $G^0$. View $\{3, 4\}$ in $G^0$ as $\{1, 2\}$ in $H_2$, view $v15$ in $G^0$ ($v\in B$) as $i34$ ($i\in D$) in $H_2$. Then $G^0$ is contained in an isomorphic copy of $H_2$, a contradiction.

So we may assume that $v23, v24\in E(G)$, then $v34\notin E(G)$. Therefore $\{v, 1, 2, 3\}, \{v, 1, 2, 4\}, \{v, 1, 3,\\5\}, \{v, 1, 4, 5\}$ form $K_4^3$. So $N^*(x, y)=\{1\}$ and $x23, x24, x35, x45\notin E(G)$ for all $x, y\in [n]\setminus\{v, 1, 2, 3, 4, 5\}$, so $B=\{v\}$. Deleting $v12, 1xv, 2xv, xy1$, adding $v34, 3xv, 4xv, xy3$ for $x, y\in E$, replacing $B_x$ by $A_x$ for $x\in E$, we obtain $G^0$. View $\{3, 4\}$ in $G^0$ as $\{1, 2\}$ in $H_1$, view $\{v15, 125\}$ in $G^0$ as $\{345, 346\}$ in $H_1$. Then $G^0$ is contained in an isomorphic copy of $H_1$, a contradiction.

{\em Case 3.2.2.} $v15, v25\in E(G)$.

If $v34\notin E(G)$, since $v\in B$, then $\{v, 1, 4, 5\}, \{v, 2, 4, 5\}$ form $K_4^3$. Since $G$ is $K_4^3\cup e$-free, then $N^*(x, y)=\emptyset$ for $x, y\in [n]\setminus\{v, 1, 2, 3, 4, 5\}$, contradicting to $G$ being dense. Therefore $v34\in E(G)$.

If $v35\notin E(G)$ (or $v45\notin E(G)$), then $\{v, 1, 3, 4\}, \{v, 2, 3, 4\}$ form $K_4^3$. Since $G$ is $K_4^3\cup e$-free, then $N^*(x, y)=\emptyset$ for $x, y\in [n]\setminus\{v, 1, 2, 3, 4, 5\}$, contradicting to $G$ being dense. Therefore $v35, v45\in E(G)$.

If $v13\notin E(G)$ (or $v14\notin E(G)$ or neither $v13$ nor $v14$ is in $E(G)$ or $v23\notin E(G)$ or $v24\notin E(G)$ or neither $v23$ nor $v24$ is in $E(G)$), then $\{v, 2, 3, 4\}, \{v, 2, 3, 5\}, \{v, 2, 4, 5\}$ form $K_4^3$. So $N^*(x, y)=\{2\}$ and $x13, x14, x15\notin E(G)$ for all $x, y\in [n]\setminus\{v, 1, 2, 3, 4, 5\}$. By {\em Case 3.2.1}, $x\notin B$ for $x\in [n]\setminus\{v, 1, 2, 3, 4, 5\}$, so $B=\{v\}$. Replacing $B_x$ by $A_x$ for $x\in E$, deleting $xy2, v15, 1xv, 2xv$, and adding $xy3, v13, 3xv, 4xv$, respectively, we obtain$G^0$. Note that $N^*(x, y)=\{3\}$ for $x, y\in E$ in $G^0$. View $\{3, 4\}$ in $G^0$ as $\{1, 2\}$ in $H_1$, view $\{v25, 125\}$ in $G^0$ as $\{345, 346\}$ in $H_1$. Then $G^0$ is contained in an isomorphic copy of $H_1$, a contradiction.

{\em Case 3.3.} There exists $v\in B$ and $v12, v15, v25\in E(G)$.

Since $v\in B$, then there are at most 2 of $\{v13, v14, v23, v24, v34, v35, v45\}$ not in $E(G)$. We claim that there are exactly 2 of those edges not in $E(G)$. Otherwise if there is at most 1 of $\{v13, v14, v23, v24, v34, \\v35, v45\}$ not in $E(G)$, without loss of generality, say at most $v13\notin E(G)$, then $\{v, 1, 4, 5\}, \{v, 2, 4, 5\}$ form $K_4^3$. So $N^*(x, y)=\emptyset$ for all $x, y\in [n]\setminus\{v, 1, 2, 3, 4, 5\}$, a contradiction. So we may assume that there are exactly 2 of $v13, v14, v23, v24, v34, v35, v45$ not in $E(G)$.

Assume that $v34\notin E(G)$. If $v13\notin E(G)$ (the case that 1 of $\{v14, v23, v24\}$ is not in $E(G)$ is similar), then $\{v, 1, 4, 5\}, \{v, 2, 3, 5\}$ form $K_4^3$. Since $G$ is $K_4^3\cup e$-free, then $N^*(x, y)=\emptyset$ for $x, y\in [n]\setminus\{v, 1, 2, 3, 4, 5\}$, contradicting to $G$ being dense. If $v35\notin E(G)$ (or $v45\notin E(G)$), then $\{v, 1, 4, 5\}, \{v, 2, 4, 5\}$ form $K_4^3$. Since $G$ is $K_4^3\cup e$-free, then $N^*(x, y)=\emptyset$ for $x, y\in [n]\setminus\{v, 1, 2, 3, 4, 5\}$, contradicting to $G$ being dense. So we may assume that $v34\in E(G)$.

Assume that 1 of $\{v13, v14, v23, v24\}$ is not in $E(G)$. Without loss of generality, let $v13\notin E(G)$.
If $v14\notin E(G)$, then $\{v, 1, 2, 5\}, \{v, 2, 3, 4\}, \{v, 2, 3, 5\}, \{v, 2, 4, 5\}$ form $K_4^3$. So $N^*(x, y)=\{2\}$ and $x13, x14, x34, x35\notin E(G)$ for all $x, y\in [n]\setminus\{v, 1, 2, 3, 4, 5\}$. Therefore $B=\{v\}$. Replacing $v12, v15$ by $v13, v14$, replacing $B_x$ by $A_x$ for $x\in E$, deleting $1xv, 2xv, xy2$ and adding $3xv, 4xv, xy3$ for $x, y\in E$, respectively. Let $G^0$ be the resulting 3-graph. Note that $N^*(x, y)=\{3\}$ for $x, y\in E$ in $G^0$. View $\{3, 4\}$ in $G^0$ as $\{1, 2\}$ in $H_1$, view $\{v25, 125\}$ in $G^0$ as $\{345, 346\}$ in $H_1$. Then $G^0$ is contained in an isomorphic copy of $H_1$, a contradiction.
If $v23\notin E(G)$ (or $v35\notin E(G)$), then $\{v, 1, 4, 5\}, \{v, 2, 4, 5\}$ form $K_4^3$. So $N^*(x, y)=\emptyset$ for all $x, y\in [n]\setminus\{v, 1, 2, 3, 4, 5\}$, a contradiction. If $v24\notin E(G)$, then $\{v, 1, 4, 5\}, \{v, 2, 3, 5\}$ form $K_4^3$. Since $G$ is $K_4^3\cup e$-free, then $N^*(x, y)=\emptyset$ for $x, y\in [n]\setminus\{v, 1, 2, 3, 4, 5\}$, contradicting to $G$ being dense.
If $v45\notin E(G)$, then $\{v, 1, 2, 4\}, \{v, 1, 2, 5\}, \{v, 2, 3, 4\}, \{v, 2, 3, 5\}$ form $K_4^3$. So $N^*(x, y)=\{2\}$ and $x14, x15, x34, x35\notin E(G)$ for all $x, y\in [n]\setminus\{v, 1, 2, 3, 4, 5\}$. Therefore $B=\{v\}$. Replacing $v12, v25$ by $v13, v45$, replacing $B_x$ by $A_x$ for $x\in E$, deleting $1xv, 2xv, xy2$ and adding $3xv, 4xv, xy3$ for $x, y\in E$, respectively. Let $G^0$ be the resulting 3-graph. Note that $N^*(x, y)=\{3\}$ for $x, y\in E$ in $G^0$. View $\{3, 4\}$ in $G^0$ as $\{1, 2\}$ in $H_1$, view $\{v15, 125\}$ in $G^0$ as $\{345, 346\}$ in $H_1$. Then $G^0$ is contained in an isomorphic copy of $H_1$, a contradiction. So it's sufficient to consider $v35, v45\notin E(G)$. However $\{v, 1, 3, 4\}, \{v, 2, 3, 4\}$ form $K_4^3$ in this situation. Since $G$ is $K_4^3\cup e$-free, then $N^*(x, y)=\emptyset$ for $x, y\in [n]\setminus\{v, 1, 2, 3, 4, 5\}$, contradicting to $G$ being dense.
\q

\subsection{$G$ does not contain two copies of $K_4^3$ sharing two vertices}
Before giving the proof of Lemma \ref{X2}, we will prove the following Lemmas.
\begin{lemma}\label{X4}
$G$ is $X_4$-free.
\end{lemma}
{\em Proof of Lemma \ref{X4}.} Assume that $G$ contains an $X_4$ with the vertex set $\{1, 2, 3, 4, 5, 6, 7, 8, 9, 10\}$, and $\{1, 2, 3, 4\}$, $\{1, 2, 5, 6\}$, $\{1, 2, 7, 8\}$, $\{1, 2, 9, 10\}$ form $K_4^3$. Since $G$ is $K_4^3\cup e$-free, then there is no edge in $V(G)\setminus\{1, 2\}$. Therefore $G$ is a subgraph of $B(2, n-2)$. By Lemma \ref{B2n}, $\lambda(G)\leq\frac{\sqrt 3}{18}$.
\q

\begin{lemma}\label{X3}
$G$ is $X_3$-free.
\end{lemma}
{\em Proof of Lemma \ref{X3}.} Assume that $G$ contains an $X_3$ with vertex set $\{1, 2, 3, 4, 5, 6, 7, 8\}$ and $\{1, 2, 3, 4\}$, $\{1, 2, 5, 6\}$, $\{1, 2, 7, 8\}$ form $K_4^3$. Denote $C=\{3, 4, 5, 6, 7, 8\}$, $D=\{x\in [n]\setminus\{1, 2, 3, 4, 5, 6, 7, 8\}: x12\in E(G)\}$ and $E=[n]\setminus(\{1, 2\}\cup C\cup D)$. Note that $x12\notin E(G)$ for $x\in E$. Since $G$ is $K_4^3\cup e$-free and $X_4$-free, then $N^*(x, y)=\{1\}$ or $\{2\}$ for $x, y\in D$, and there is no edge between $C$ and $D\cup E$. And if there is an edge $e_1$ in $C$, then $|e_1\cap\{3, 4\}|$=$|e_1\cap\{5, 6\}|$=$|e_1\cap\{7, 8\}|=1$, i.e. $G[C]$ is a 3-partite 3-graph and $\lambda(G[C])\leq\frac{1}{27}$. Set
$$x_1=a, \ x_2=b, \ \sum_{v\in C}x_v=c, \ \sum_{v\in D}x_v=d \ and \ \sum_{v\in E}x_v=e.$$
Without loss of generality, assume that $a\geq b$, then replacing $xy2$ by $xy1$ for all $x, y\in D$ does not decrease the Lagrangian. So
\begin{eqnarray*}
\lambda(G)&\leq&ab(c+d)+a(\bigg(\frac{c}{6}\bigg)^2{6\choose 2}+cd+ce+de+\frac{d^2}{2}+\frac{e^2}{2})+b(\bigg(\frac{c}{6}\bigg)^2{6\choose 2}+cd+ce+de+\frac{e^2}{2})+\frac{c^3}{27}\\
&\leq&ab(c+d)+a\bigg(\frac{5c^2}{12}+c(d+e)+\frac{(d+e)^2}{2}\bigg)+b\bigg(\frac{5c^2}{12}+c(d+e)+\frac{(d+e)^2}{2}\bigg)+\frac{c^3}{27}.
\end{eqnarray*}
Let $\delta=d+e$. Then
$$\lambda(G)\leq\lambda(a, b, c, \delta)=ab(c+\delta)+a(\frac{5c^2}{12}+c\delta+\frac{\delta^2}{2})+b(\frac{5c^2}{12}+c\delta+\frac{\delta^2}{2})+\frac{c^3}{27}\triangleq\lambda$$
subject to
\begin{equation}
\left\{
             \begin{array}{lr}
             a+b+c+\delta=1, &  \\
             a\geq 0, \ b\geq 0, \ c\geq 0 \ and \ \delta\geq 0.&
             \end{array}
\right.
\end{equation}
For simplicity of the notation, we assume that $\lambda$ reaches the maximum at $(a, b, c, \delta)$.

If $c=0$, then $\lambda\leq\frac{(a+b)^2[1-(a+b)]}{4}+\frac{(a+b)[1-(a+b)]^2}{2}$. By Fact \ref{aaa}, $\lambda\leq\frac{\sqrt 3}{18}.$

If $\delta=0$, then
\begin{eqnarray*}
\lambda&=&abc+\frac{5(a+b)c^2}{12}+\frac{c^3}{27}\\
&\leq&\frac{(a+b)^2c}{4}+\frac{5(a+b)c^2}{12}+\frac{c^3}{27}\\
&=&\frac{(1-c)^2c}{4}+\frac{5(1-c)c^2}{12}+\frac{c^3}{27}=f(c).\\
f'(c)&=&\frac{(1-c)^2-2c(1-c)}{4}+\frac{10c(1-c)-5c^2}{12}+\frac{c^2}{9}\\
&=&\frac{-14c^2-6c+9}{36}.
\end{eqnarray*}
So $f(c)$ is increasing in $[0, \frac{3(\sqrt {15} -1)}{14}]$, then $f_{max}=f(\frac{3(\sqrt {15} -1)}{14})<0.0921.$

If $c, \delta>0$, then by Theorem \ref{KKT}, we have $\frac{\partial\lambda}{\partial c}=\frac{\partial\lambda}{\partial \delta},$ solving it, we obtain that  $\frac{ac+bc}{6}=\frac{c^2}{9}.$ Therefore $c=\frac{3(a+b)}{2}$. So
\begin{eqnarray*}
\lambda&=&ab(1-a-b)+\frac{(a+b)(c+\delta)^2}{2}-\frac{c^2(a+b)}{12}+\frac{c^3}{27}\\
&\leq&\frac{(a+b)^2[1-(a+b)]}{4}+\frac{(a+b)[1-(a+b)]^2}{2}+\frac{c^3}{27}-\frac{c^3}{18}\\
&<&\frac{(a+b)^2[1-(a+b)]}{4}+\frac{(a+b)[1-(a+b)]^2}{2}\\
&\leq&\frac{\sqrt 3}{18},
\end{eqnarray*}
the last inequality follows from Fact \ref{aaa}.
\q

{\em Proof of Lemma \ref{X2}.} Assume that $G$ contains an $X_2$ with vertex set $\{1, 2, 3, 4, 5, 6\}$ and $\{1, 2, 3, 4\}$, $\{1, 2, 5, 6\}$ form $K_4^3$. We prove the following claim first.
\begin{claim}\label{X2''}
$\lambda(G[\{1, 2, 3, 4, 5, 6\}])\leq \frac{2}{25}.$
\end{claim}
{\em Proof of Claim \ref{X2''}.} Denote $A=\{1, 2, 3, 4, 5, 6\}$. By Lemma \ref{K53-}, $G$ is $K_5^{3-}$-free, then $e(G[A'])\leq 8$ for all $A'\subseteq A$ and $|A'|=5$. By double counting, we have $${6-3\choose2}e(G[A])\leq{6\choose 5}\times8.$$
Therefore $e(G[A])\leq16={6-1\choose3}+{6-2\choose2}<{6\choose 3}-3=17$. By Lemma \ref{TPP} and Lemma \ref{T}, $\lambda(G[A])\leq\lambda(K_5^3)=\frac{2}{25}.$
\q

Denote $C=\{3, 4, 5, 6\}$, $D=\{x\in V(G)\setminus\{1, 2, 3, 4, 5, 6\}: x12\in E(G)\}$ and $E=V(G)\setminus(\{1, 2\}\cup C\cup D)$. Note that $x12\notin E(G)$ for $x\in E$. Since $G$ is $K_4^3\cup e$-free, then $x34, x56\notin E(G)$ and only possibly, $x35, x36, x45, x46\in E(G)$ for $x\in D\cup E$. Let
$$x_1=a, x_2=b, \sum_{v\in C}x_v=c, \sum_{v\in D}x_v=d, \sum_{v\in E}x_v=e.$$
Without loss of generality, assume that $a\geq b$. Note that the contribution of the edges between $C$ and $D\cup E$ to $\lambda(G)$ is at most $(d+e)(x_3+x_4)(x_5+x_6)\leq\frac{(d+e)c^2}{4}$. Since $G$ is $X_3$-free, then $N^*(x, y)=\{1\}$ or $\{2\}$ for $x, y\in D$. If $xy2\in E(G)$, we delete $xy2$ and add $xy1$, this does not reduce the Lagrangian. Hence
\begin{eqnarray}
\lambda(G)&\leq&abd+a(cd+ce+de+\frac{d^2}{2}+\frac{e^2}{2})+b(cd+ce+de+\frac{e^2}{2})+\frac{2(a+b+c)^3}{25}+\frac{(d+e)c^2}{4} \label{eq1}\\
&=&\lambda(a, b, c, d, e)=\lambda \notag
\end{eqnarray}
under the constraints $a+b+c+d+e=1$, $a, b, c, d, e\geq 0$.

To simplify the notation, we assume that $\lambda$ reaches the maximum at $(a, b, c, d, e)$, Note that $a\geq b$.

\begin{claim}\label{b0}
$\lambda(a, 0, c, d, e)\leq\frac{\sqrt 3}{18}$.
\end{claim}
The proof of Claim \ref{b0} is given in Appendix.

Let us continue the proof of Lemma \ref{X2}. We have shown that $a\geq b>0$ (Claim \ref{b0}). If $d=0$, substitute it into (\ref{eq1}), then
$$\lambda(a, b, c, 0, e)=a(ce+\frac{e^2}{2})+b(ce+\frac{e^2}{2})+\frac{2(a+b+c)^3}{25}+\frac{ec^2}{4}.$$
So $\lambda(a+b, 0, c, 0, e)$ also gets the maximum value, a contradiction to $a\geq b>0$ when $\lambda$ gets the maximum. So $a, b, d>0$. By Theorem \ref{KKT}, $\frac{\partial\lambda}{\partial a}=\frac{\partial\lambda}{\partial b}$, combining with
\begin{eqnarray*}
\frac{\partial\lambda}{\partial a}&=&bd+cd+ce+de+\frac{d^2}{2}+\frac{e^2}{2}+\frac{6(a+b+c)^2}{25},\\
\frac{\partial\lambda}{\partial b}&=&ad+cd+ce+de+\frac{e^2}{2}+\frac{6(a+b+c)^2}{25},
\end{eqnarray*}
we get
\begin{eqnarray}
a=b+\frac{d}{2}. \label{eqa}
\end{eqnarray}

\begin{claim}\label{e0}
$\lambda(a, b, c, d, 0)\leq\frac{\sqrt 3}{18}$.
\end{claim}
The proof of Claim \ref{e0} is given in Appendix.

Let us continue the proof of Lemma \ref{X2}. We have shown that $a, b, d, e>0$. By Theorem \ref{KKT}, we have $\frac{\partial\lambda}{\partial d}=\frac{\partial\lambda}{\partial e}$. Since
\begin{eqnarray*}
\frac{\partial\lambda}{\partial d}&=&ab+ac+ae+ad+bc+be+\frac{c^2}{4},\\
\frac{\partial\lambda}{\partial e}&=&ac+ad+ae+bc+bd+be+\frac{c^2}{4},
\end{eqnarray*}
then $a=d$. Recall  (\ref{eqa}), so
\begin{eqnarray}
a=d=2b. \label{abd}
\end{eqnarray}

\begin{claim}\label{c0}
$\lambda(a, b, 0, d, e)\leq\frac{\sqrt 3}{18}$.
\end{claim}
{\em Proof of Claim \ref{c0}.}
Assume $c=0$. By (\ref{abd}), then $e=1-a-b-c-d=1-5b$. Substituting it into (\ref{eq1}), then
\begin{eqnarray*}
\lambda&=&abd+a(de+\frac{d^2}{2}+\frac{e^2}{2})+b(de+\frac{e^2}{2})+\frac{2(a+b)^3}{25}\\
&=&4b^3+b\bigg(6b(1-5b)+4b^2+\frac{3(1-5b)^2}{2}\bigg)+\frac{2(3b)^3}{25}\\
&=&\frac{883b^3}{50}-9b^2+\frac{3b}{2}=f(b),\\
f'(b)&=&\frac{2649b^2}{50}-18b+\frac{3}{2}.
\end{eqnarray*}
So $f(b)$ is increasing in $[0, \frac{450-5\sqrt{153}}{2649}]$ or $[\frac{450+5\sqrt{153}}{2649}, 1]$. Note that $a+b+d+e=1$, then $5b<1$, so $\lambda\leq f(\frac{450-5\sqrt{153}}{2649})\leq 0.083$.
\q

Let us continue the proof of Lemma \ref{X2}. We have shown that $a=d=2b$ and $e=1-5b-c$ and $a, b, c, d, e>0$. Substituting them into (\ref{eq1}), we have
\begin{eqnarray*}
\lambda&=&4b^3+b(6bc+3ce+6be+4b^2+\frac{3e^2}{2})+\frac{2(3b+c)^3}{25}+\frac{(2b+e)c^2}{4}\\
&=&\frac{31b^3}{2}-\frac{3bc^2}{2}-9b^2+\frac{3b}{2}+\frac{2(3b+c)^3}{25}+\frac{(1-3b-c)c^2}{4}\\
&\leq&16b^3-\frac{3bc^2}{2}-9b^2+\frac{3b}{2}+\frac{(3b+c)^3}{12}+\frac{(1-3b-c)c^2}{4}=\lambda_0(b, c)
\end{eqnarray*}
under the constraints $5b+c\leq1, b, c\ge 0.$ Now we estimate the optimum value of $\lambda_0$. For simplicity of the notation, let $\lambda_0$ reach the maximum value at $(b, c)$.

If $b=0$, then $\lambda_0=\frac{c^3}{12}+\frac{(1-c)c^2}{4}\leq\frac{3c^2-2c^3}{12}=f(c)$. Since $f'(c)=\frac{c(1-c)}{2}$, then $\lambda$ is increasing in $[0, 1]$. Therefore $\lambda_0\leq f(1)=\frac{1}{12}<0.09.$

If $c=0$, then $\lambda_0=16b^3-9b^2+\frac{3b}{2}+\frac{(3b)^3}{12}=\frac{73b^3}{4}-9b^2+\frac{3b}{2}\leq19b^3-9b^2+\frac{3b}{2}=f(b)$. Since $f'(b)=57b^2-18b+\frac{3}{2}>0$, so $f(b)$ is increasing in $[0, 1]$. Therefore $\lambda_0\leq f(\frac{1}{5})=0.092.$

If $5b+c=1$, then
\begin{eqnarray*}
\lambda_0&=&16b^3-\frac{3b(1-5b)^2}{2}-9b^2+\frac{3b}{2}+\frac{(1-2b)^3}{12}+\frac{b(1-5b)^2}{2},\\
\lambda_0'&=&-29b^2+4b.
\end{eqnarray*}
Then $\lambda_0\leq \lambda_0(4/29)\leq 0.0961<\frac{\sqrt 3}{18}$.

Therefore we may assume that $\lambda_0$ gets maximum when $b, c>0$ and $5b+c<1$. By Theorem \ref{KKT},
\begin{eqnarray*}
\frac{\partial\lambda_0}{\partial b}&=&48b^2-\frac{9c^2}{4}-18b+\frac{3}{2}+\frac{3(3b+c)^2}{4}=0\\
\frac{\partial\lambda_0}{\partial c}&=&-3bc+\frac{(3b+c)^2}{4}+\frac{2c-6bc-3c^2}{4}=0.
\end{eqnarray*}
Equivalently, $\frac{\partial\lambda_0}{\partial b}-3\times\frac{\partial\lambda_0}{\partial c}=0$ and $\frac{\partial\lambda_0}{\partial c}=0$. Solving these two equations, we have
$$c=\frac{1-12b+32b^2}{1-9b}=\frac{(1-8b)(1-4b)}{1-9b}\quad {\rm and}\quad 9b^2-12bc+2c-2c^2=0.$$
Recall that $c>0$ and $b<\frac{1}{5}$. So $0<b<\frac{1}{9}$ or $\frac{1}{8}\leq b\leq \frac{1}{5}$. Combining the above equations, we have $b(2137b^3-882b^2+125b-6)=0$. Let $f(b)=2137b^3-882b^2+125b-6$. Since $f'(b)=6411b^2-1764b+125>0$, then $f(b)$ is increasing in $[0, \frac{1}{4}]$. However $f(0), f(\frac{1}{9})<0$ and $f(\frac{1}{8})>0$, so $f(b)=0$ has no solution in $0<b<\frac{1}{9}$ and $\frac{1}{8}\leq b\leq \frac{1}{5}$, a contradiction. This completes the proof of Lemma \ref{X2}.
\q

\subsection{$G$ contains at least two copies of $K_4^3$}
In this section, we give the proof of Lemma \ref{2K4}.

{\em Proof of Lemma \ref{2K4}.} Recall that $\{1, 2, 3, 4\}$ forms a $K_4^3$. Assume that $G$ contains no other $K_4^3$, in other words, $v$ does not belong to any $K_4^3$ for any $v\in [n]\setminus\{1, 2, 3, 4\}$. We claim that $|G_v\cap\{12, 13, 14, 23, 24, 34\}|\leq4$. Since otherwise $G_v[\{1, 2, 3, 4\}]$ contains a triangle and $v$ is contained in the $K_4^3$ formed by $v$ and the vertices in this triangle.
Since $G$ is $K_4^3\cup e$-free, then $G_v$ does not contain an edge in $[n]\setminus\{1, 2, 3, 4\}$. By Claim \ref{clique}, $\omega(G_v)\geq 3$, so the maximum clique of $G_v$ contains at least 2 vertices in $\{1, 2, 3, 4\}$, therefore $|G_v\cap\{12, 13, 14, 23, 24, 34\}|\geq1$.
Let $$A=\{1, 2, 3, 4\},$$ $$A_1=\{v\in [n]\setminus A: |G_v\cap\{12, 13, 14, 23, 24, 34\}|=1\},$$
$$A_2=\{v\in [n]\setminus A: |G_v\cap\{12, 13, 14, 23, 24, 34\}|=2\},$$
$$A_3=\{v\in [n]\setminus A: |G_v\cap\{12, 13, 14, 23, 24, 34\}|=3\},$$
$$A_4=\{v\in [n]\setminus A: |G_v\cap\{12, 13, 14, 23, 24, 34\}|=4\}.$$
Without loss of generality, let's assume that $x_1\geq x_2\geq x_3\geq x_4$. Then $x_1x_2\geq x_1x_3\geq x_2x_3, x_1x_4\geq x_2x_4\geq x_3x_4.$ We aim to give an upper bound of $\lambda(G, \vec x)$, therefore we can assume that $v12, v13\in E(G)$ for $v\in A_2$, $v12\in E(G)$ for $v\in A_1$.
Set $$x_1=a, x_2=b, x_3=c, x_4=d, \sum_{v\in A_1}x_v=h, \sum_{v\in A_2}x_v=g, \sum_{v\in A_3}x_v=f, \sum_{v\in A_4}x_v=e.$$
Since $G$ doesn't contain two copies of $K_4^3$, then the deletion of any 1 of $\{123, 124, 134, 234\}$ of $E(G)$ makes $G$ $K_4^3$-free. So $abc, abd, acd, bcd>0.00264$ since otherwise $\lambda(G)\leq 0.00264+\frac{0.5615}{6}$(in view of Lemma \ref{K43})$\leq\frac{\sqrt 3}{18}.$ So $(a+b)cd>2\times 0.00264$. By Claim \ref{2}, $a+b\leq\frac{3-\sqrt 3}{3}$, then $cd\geq\frac{6\times 0.00264}{3-\sqrt 3}$. Therefore $c+d\geq2\sqrt{cd}>0.22354.$ If $d<0.11177$, then $b+c\geq2\sqrt{bc}>0.307.$
%By Claim \ref{clique}, then $a\leq1-\sqrt{\frac{\sqrt 3}{3}}\leq 0.2402.$

 To complete the proof, we show the following three claims in Appendix.
\begin{claim}\label{A4}
$\lambda(G[A\cup A_4], \vec x)\leq 0.0789(a+b+c+d+e)^3.$
\end{claim}

\begin{claim}\label{A3}
$\lambda(G[A\cup A_4\cup A_3], \vec x)\leq 0.092(a+b+c+d+e+f)^3.$
\end{claim}

\begin{claim}\label{A2}
$\lambda(G)\leq \frac{\sqrt 3}{18}.$
\end{claim}

\subsection{$G$ does not contain two copies of $K_4^3$ sharing three vertices}
{\em Proof of Lemma \ref{Y2}.} Assume that $G$ contains an $Y_2$ with the vertex set $\{1, 2, 3, 4, 5\}$, where $\{1, 2, 3, 4\}$ and $\{1, 2, 3, 5\}$ form $K_4^3$. Since $G$ is $K_4^3\cup e$-free, then any two $K_4^3$ in $G$ must intersect 2 or 3 vertices. Since $G$ is $X_2$-free (Lemma \ref{X2}), then any two $K_4^3$ in $G$ must intersect 3 vertices. Therefore $G-\{3\}$ cannot contain a $K_4^3$ since it cannot intersect with $\{1, 2, 3, 4\}$ and $\{1, 2, 3, 5\}$ at three vertices in the same time. Let $x_1=a$, $x_2=b$, $x_3=c$, $x_4+x_5=d$ and assume that $a\geq b\geq c$.  By Claim \ref{3}, $c>0.08$.
%Since $G$ is $K_4^3\cup\{e\}$-free, $K_5^3$-free and $X_2$-free, then $E(G)-\{123\}$ is $K_4^3$-free. Therefore $$\lambda(G)\leq abc+\frac{0.5616}{6}.$$
%If $\lambda(G)\leq\frac{\sqrt 3}{18}$, then we are done. Otherwise $\lambda(G)\geq\frac{\sqrt 3}{18}$, then $abc>\frac{\sqrt 3}{18}-\frac{0.5615}{6}>0.0026417$. By Claim \ref{2}, then $a+b\leq\frac{3-\sqrt 3}{3}$. Therefore $ab\leq(\frac{3-\sqrt 3}{6})^2$, then we can always assume that $c>0.059$.
Let $D=\{4, 5\}$, and
$$E_0=\{v: v\in [n]\setminus\{1, 2, 3, 4, 5\} \ and \ |G_v\cap\{12, 13, 23\}|=0\},$$
$$E_1=\{v: v\in [n]\setminus\{1, 2, 3, 4, 5\} \ and \ |G_v\cap\{12, 13, 23\}|=1\},$$
$$E_2=\{v: v\in [n]\setminus\{1, 2, 3, 4, 5\} \ and \ |G_v\cap\{12, 13, 23\}|\geq 2\}.$$
Set $\sum_{v\in E_0}x_v=g$, $\sum_{v\in E_1}x_v=f$ and $\sum_{v\in E_2}x_v=e$. Since $G$ is $K_4^3\cup e$-free, then for $x, y\in [n]\setminus\{1, 2, 3, 4, 5\}$, $N^*(x, y)\subseteq\{1, 2, 3\}$. If $x\in D$ and $y\in E_2$, then $N^*(x, y)\subseteq\{1, 2, 3, 4, 5\}$. We claim that $|N^*(x, y)\cap\{1, 2, 3\}|\leq 2$ for $x, y\in E_2$ or $x\in D$ and $y\in E_2$, otherwise there exists two vertices $z, w\in N^*(x, y)\cap\{1, 2, 3\}$ such that $zw\in G_x\cap G_y$. So $\{x, y, z, w\}$ forms a $K_4^3$, which forms an $X_2$ with $G[\{1, 2, 3, 4\}]$, a contradiction. Therefore we may let $N^*(x, y)\cap\{1, 2, 3\}\subseteq\{1, 2\}$ for those $x, y$ with Lagrangian non-decreasing. Since $G$ is $K_4^3\cup e$-free, then all edges in $G[V-\{1, 2, 3\}]$ must contain $\{4, 5\}$. Since $G$ is $K_5^3$-free, then at least one of $\{145, 245, 345\}$ is not in $E(G)$. (Indeed, $G$ is $K_5^{3-}$-free, at least two of $\{145, 245, 345\}$ are not in $E(G)$. But it seems to be easier to estimate the Lagrangian below if we relax it to be one.) We may assume that 345 is not in $E(G)$. Therefore
\begin{eqnarray*}
\lambda(G)&\leq&ab(c+d+e+f)+(ac+bc)(d+e)+(a+b)(de+\frac{e^2}{2})+(a+b+c)(df+dg+eg+fg+ef\\
&+&\frac{f^2}{2}+\frac{g^2}{2})+\frac{d^2(a+b+e+f+g)}{4}\\
%ab(c+d+e+f)+(ac+bc)(d+e)+(a+b)(de+df+\frac{e^2}{2})+(a+b+c)(dg+eg+fg+ef\\
%&+&\frac{f^2}{2}+\frac{g^2}{2})+\frac{d^2(a+e+f+g)}{4}\\
&=&ab(c+d+e+f)+(ac+bc)(d+e)+(a+b)(de+\frac{e^2}{2})+(a+b+c)\bigg(d(g+f)+e(g+f)\\
&+&\frac{(f+g)^2}{2}\bigg)+\frac{d^2(a+b+e+f+g)}{4}\\
&=&\lambda(a, b, c, d, e, f, g)=\lambda.
\end{eqnarray*}
Note that $\lambda(a, b, c, d, e, f, g)\leq \lambda(\frac{a+b}{2}, \frac{a+b}{2}, c, d, e, f+g, 0)$, then we may  assume that $g=0$ and $a=b$. Then let $\alpha=a+b$, so
\begin{eqnarray*}
\lambda&=&\frac{\alpha^2(c+d+e+f)}{4}+\alpha c(d+e)+\alpha(de+\frac{e^2}{2})+(\alpha+c)(df+ef+\frac{f^2}{2})+\frac{d^2(\alpha+e+f)}{4}
\end{eqnarray*}
subject to
\begin{equation}
\left\{
             \begin{array}{lr}
             \alpha+c+d+e+f=1, &  \\
             c\geq0.08, \alpha, d, e, f\ge 0. & \label{eq6}
             \end{array}
\right.
\end{equation}
Note that
\begin{eqnarray*}
\lambda&=&\frac{\alpha^2(c+d+e+f)}{4}+\alpha c(d+e)+\alpha(\frac{d^2}{2}+de+\frac{e^2}{2})+(\alpha+c)(df+ef+\frac{f^2}{2})+\frac{d^2(e+f-\alpha)}{4}.
\end{eqnarray*}
If $e+f-\alpha<\frac{d}{2}$, then $\lambda\leq\lambda(a, c, d-\epsilon, e+\epsilon, f)$ for $\epsilon>0$ small enough, a contradiction. So we may  assume that $e+f\geq\alpha+\frac{d}{2}$ or $d=0$.

\begin{claim}\label{2a0}
$\lambda(0, c, d, e, f)\leq\frac{\sqrt 3}{18}$.
\end{claim}
{\em Proof of Claim \ref{2a0}.} If $\alpha=0$, then $c+d+e+f=1$ and
\begin{eqnarray*}
\lambda&=&c(df+ef+\frac{f^2}{2})+\frac{d^2(e+f)}{4}\\
&\leq&c\frac{(d+e+f)^2}{2}+\bigg(\frac{\frac{d}{2}+\frac{d}{2}+e+f}{3}\bigg)^3\\
&=&\frac{c(1-c)^2}{2}+\frac{(1-c)^3}{27}=\frac{25c^3-48c^2+21c+2}{54}=f(c),\\
f'(c)&=&\frac{25c^2-32c+7}{18}=\frac{(25c-7)(c-1)}{18}.
\end{eqnarray*}
Therefore $\lambda\leq f(\frac{7}{25})=0.0864.$
\q

Let's continue the proof of Lemma \ref{Y2}. We claim that if $c>0.08$, then $\alpha\geq c$.  If $\alpha<c$, then $\alpha\leq\frac{1}{2}$, and for $0<\epsilon<\frac{a}{3}$,
\begin{eqnarray*}
\lambda(\alpha+\epsilon, c-\epsilon, d, e, f)-\lambda(\alpha, c, d, e, f)&\geq&\frac{(\alpha+\epsilon)^2(c-\epsilon+d+e+f)-(\alpha)^2(c+d+e+f)}{4}\\
&>&\epsilon(2\alpha-3\alpha^2-3\epsilon\alpha)\\
&>&\epsilon(2\alpha-4\alpha^2)\geq0,
\end{eqnarray*}
a contradiction.

\begin{claim}\label{2f0}
$\lambda(\alpha, c, d, e, 0)\leq\frac{\sqrt 3}{18}$.
\end{claim}
The proof of Claim \ref{2f0} is given in Appendix.

\begin{claim}\label{2e0}
$\lambda(\alpha, c, d, 0, f)\leq\frac{\sqrt 3}{18}$.
\end{claim}
The proof of Claim \ref{2e0} is given in Appendix.

By the above claims, we may assume that $\alpha, e, f>0$, so by Theorem \ref{KKT}, we have $\frac{\partial\lambda}{\partial e}=\frac{\partial\lambda}{\partial f}$. In view of (\ref{eq6}),
\begin{eqnarray*}
%\frac{\partial \lambda}{\partial a}&=&\frac{a(c+d+e+f)}{2}+cd+ce+de+\frac{e^2}{2}+df+ef+\frac{f^2}{2}+\frac{d^2}{4}\\
%\frac{\partial \lambda}{\partial c}&=&\frac{a^2}{4}+ad+ae+df+ef+\frac{f^2}{2}\\
%\frac{\partial \lambda}{\partial d}&=&\frac{a^2}{4}+ac+ae+(a+c)f+\frac{d(a+e+f)}{2}\\
\frac{\partial \lambda}{\partial e}&=&\frac{\alpha^2}{4}+\alpha c+\alpha d+\alpha e+(\alpha+c)f+\frac{d^2}{4}\\
\frac{\partial \lambda}{\partial f}&=&\frac{\alpha^2}{4}+(\alpha+c)(d+e+f)+\frac{d^2}{4}.
\end{eqnarray*}
So $\alpha=d+e.$

\begin{claim}\label{2d0}
$\lambda(\alpha, c, 0, e, f)\leq\frac{\sqrt 3}{18}$.
\end{claim}
The proof of Claim \ref{2d0} is given in Appendix.

Let's continue the proof of Lemma \ref{Y2}. The above claims indicate that  we may assume $d, e>0$, then by Theorem \ref{KKT}, we have $\frac{\partial\lambda}{\partial d}=\frac{\partial\lambda}{\partial e}$. Since
\begin{eqnarray*}
\frac{\partial \lambda}{\partial d}&=&\frac{\alpha^2}{4}+\alpha c+\alpha e+(\alpha+c)f+\frac{d(\alpha+e+f)}{2}\\
\frac{\partial \lambda}{\partial e}&=&\frac{\alpha^2}{4}+\alpha c+\alpha d+\alpha e+(\alpha+c)f+\frac{d^2}{4},
\end{eqnarray*}
then $\alpha+\frac{d}{2}=e+f$. Since $\alpha=d+e$, then $f=\frac{3d}{2}$. Note that $\alpha+c+d+e+f=1$, then $f=1-2\alpha-c$ and $d=\frac{2(1-2\alpha-c)}{3}>0$ and $e=\frac{7\alpha+2c-2}{3}$. Substituting these into (\ref{eq6}), we have
\begin{eqnarray*}
\lambda&=&\frac{\alpha^2(1-\alpha)}{4}+\alpha^2c+\alpha((\frac{2(1-2\alpha-c)}{3})(\frac{7\alpha+2c-2}{3})+\frac{(\frac{7\alpha+2c-2}{3})^2}{2})\\
&+&(\alpha+c)(\alpha(1-2\alpha-c)+\frac{(1-2\alpha-c)^2}{2})+\frac{(1-2\alpha-c)^2}{9}(1-c-\frac{2(1-2\alpha-c)}{3})\\
&=&\frac{-5\alpha^3}{108}+\frac{14\alpha^2c}{9}-\frac{11\alpha^2}{36}+\frac{23\alpha c^2}{18}-\frac{14\alpha c}{9}+\frac{5\alpha}{18}+\frac{25c^3}{54}-\frac{8c^2}{9}+\frac{7c}{18}+\frac{1}{27}\\
&=&\lambda(\alpha, c).
\end{eqnarray*}
If $c=0.08$, then
\begin{eqnarray*}
\lambda&=&\frac{-5\alpha^3}{108}-\frac{1.63\alpha^2}{9}+\frac{1.4536\alpha}{9}+\frac{1.6928}{27},\\
\lambda'&=&\frac{-5\alpha^2}{36}-\frac{3.26\alpha}{9}+\frac{1.4536}{9}.
\end{eqnarray*}
Note that $\lambda$ is increasing in $[0, \frac{3\sqrt{4971}-163}{125}]$, then $\lambda\leq 0.096.$ So we may assume that $c>0.08$. Recall that $\alpha\geq c$ when $c>0.08$ and $e=\frac{7\alpha+2c-2}{3}, f=1-2\alpha-c>0$, so $\alpha>\frac{2}{9}$ and $2\alpha+c<1$. So we can maximize $\lambda(\alpha, c)$ subject to
\begin{equation}
\left\{
             \begin{array}{lr}
             2\alpha+c\leq1, &  \\
             \alpha\geq\frac{2}{9}, & \\
             c\geq 0.08. &
             \end{array}
\right.
\end{equation}
Consider
\begin{eqnarray*}
&\lambda&(\alpha+{c-0.08\over 2}, 0.08)-\lambda(\alpha, c)\\
&=&{(2-25c)(32500\alpha^2+26250\alpha c-25500\alpha+9375c^2-16150c+4728)\over 500000}=\lambda_0,\\
&\lambda_0'|_{\alpha}&=\frac{(2-25c)(260\alpha+105c-102)}{2000}.
\end{eqnarray*}
If $260\alpha+105c-102>0$, then $\lambda_0'|_{\alpha}<0$. So $\lambda_0\geq\lambda_0(\frac{1-c}{2}, c)$. Therefore
\begin{eqnarray*}
\lambda_0&\geq&\frac{-7c^3}{32}+\frac{11c^2}{32}-\frac{c}{32}+\frac{103}{250000}=\lambda_1,\\
\lambda_1'&=&\frac{-21c^2}{32}+\frac{11c}{16}-\frac{1}{32}.
\end{eqnarray*}
Note that $\lambda_1$ is increasing in $[0.08, 1]$, so $\lambda_1\geq\lambda_1(0.08)=0$. So $\lambda$ gets maximum when $c=0.08$, a contradiction.

If $260\alpha+105c-102<0$, then $\lambda_0'|_a>0$. So $\lambda_0\geq\lambda_0(\frac{2}{9}, c)$. Therefore
\begin{eqnarray*}
\lambda_0&\geq&{(2-25c)(9375c^2-{30950c\over3}+{53968 \over 81})\over500000}=\lambda_2\\
\lambda_2'&=&\frac{83c}{75}-\frac{45c^2}{32}-{6041\over 81000}.
\end{eqnarray*}
Note that $\lambda_2$ is increasing in $[0.08, 0.7]$, then $\lambda_2\geq\lambda_2(0.08)=0$. So $\lambda$ gets maximum when $c=0.08$, a contradiction. The proof of Lemma \ref{Y2} is completed.
\q

\section{Remark}
Let $\Lambda_{t}^{(r)}=\{\pi_{\lambda}(\cal F): \cal F {\rm \ is \ a \ family \ of } \ r{\rm-uniform \ graphs \ and} \ |\cal F|\le t\}.$
Proposition \ref{relationlt} (Proposition \ref{relationlt} can be generalized to a family of $r$-graphs) implies that $\Lambda_{t}^{(r)} \subseteq \Pi_{t}^{(r)}$.
\begin{ques}
Is $\Lambda_{t}^{(r)}$ the same as  $\Pi_{t}^{(r)}$?
\end{ques}

Let us propose the following conjecture implying that there exists an $r$-graph whose Tur\'an density is an irrational number.
\begin{con}
If $c\cdot{r! \over r^r}$ is in $\Pi_{1}^{(r)}$ for $r\ge 2$, then  $c\cdot{p! \over p^p}$ is in $\Pi_{1}^{(p)}$ for  $p\ge r$.
\end{con}
%In \cite{Pikhurko2}, Pikhurko showed that the set $\Pi_{\infty}^{(r)}$ has cardinality of the continuum and is closed for $r\ge 3$. This implies that an open problem related to non-jumping numbers in \cite{fprt} is equivalent to the following question.
%\begin{ques}(Frankl, Peng, R\"odl, and Talbot \cite{fprt}; Pikhurko\cite{Pikhurko2} )
%Let $r\ge 3$. Is there $a<1$ such that $[a, 1]\subset \Pi_{\infty}^{(r)}$?
%\end{ques}

%\bigskip

\section{Appendix}\label{appendix}
We  give theoretical proofs for Claims \ref{b0}, \ref{e0}, \ref{A4}-\ref{A2} and  \ref{2f0}-\ref{2d0} in this section, we have also used Lingo to run the optimization problems. The outcome by Lingo is consistent with the expected optimum values. We can provide the programming upon request.
\subsection{\em Proof of Claim \ref{b0}}
Substitute $b=0$ into (\ref{eq1}), then
\begin{eqnarray}
\lambda&=&a(cd+ce+de+\frac{d^2}{2}+\frac{e^2}{2})+\frac{2(a+c)^3}{25}+\frac{(d+e)c^2}{4},\notag\\
&=&a[c(d+e)+\frac{(d+e)^2}{2}]+\frac{2(a+c)^3}{25}+\frac{(d+e)c^2}{4}, \label{412}\\
&=&a[c\delta+\frac{\delta^2}{2}]+\frac{2(a+c)^3}{25}+\frac{\delta c^2}{4},\notag
\end{eqnarray}
where $\delta=d+e$, then $a+c+\delta=1$.

If $a=0$, then
\begin{eqnarray*}
\lambda&=&\frac{2c^3}{25}+\frac{\delta c^2}{4}\leq\frac{c^3}{12}+\frac{(1-c)c^2}{4}=\frac{3c^2-2c^3}{12}=f(c)\\
f'(c)&=&\frac{c-c^2}{2}=\frac{c(1-c)}{2}.
\end{eqnarray*}
So $f(c)$ is increasing in $[0, 1]$, then $\lambda\leq f(1)=\frac{1}{12}<\frac{\sqrt 3}{18}.$

If $\delta=0$, then $$\lambda=\frac{2(a+c)^3}{25}\leq \frac{2}{25}.$$

If $c=0$, then
\begin{eqnarray*}
\lambda&=&\frac{a\delta^2}{2}+\frac{2a^3}{25}\leq\frac{a(1-a)^2}{2}+\frac{a^3}{12}=\frac{7a^3-12a^2+6a}{12}=f(a)\\
f'(a)&=&\frac{7a^2-8a+2}{4}.
\end{eqnarray*}
So $f(a)$ is increasing in $[0, \frac{4-\sqrt 2}{7}]$ or $[\frac{4+\sqrt 2}{7}, 1]$, then $\lambda\leq max\{f(1), f(\frac{4-\sqrt 2}{7})\}.$ Note that $f(\frac{4-\sqrt 2}{7})\leq 0.078$, then $\lambda\leq f(1)=\frac{1}{12}\leq\frac{\sqrt 3}{18}.$

Therefore we may assume that $a, c, \delta>0$. Substituting $c=1-a-\delta$ into (\ref{412}), then
\begin{eqnarray*}
\lambda&=&a[(1-a-\delta)\delta+\frac{\delta^2}{2}]+\frac{2(1-\delta)^3}{25}+\frac{\delta(1-a-\delta)^2}{4}
\end{eqnarray*}
gets its maximum inside  interior points. By Theorem \ref{KKT},
\begin{eqnarray*}
\frac{\partial\lambda}{\partial a}&=&\frac{\delta}{2}-\frac{3a\delta}{2}=0,\\
\frac{\partial\lambda}{\partial \delta}&=&\frac{(1-a-\delta)(1+3a-3\delta)}{4}-\frac{6(1-\delta)^2}{25}=0.
\end{eqnarray*}
Then $a=\frac{1}{3}$ and $\delta=\frac{26-10\sqrt 2}{51}$, and $\lambda< 0.09.$
\q

\subsection{\em Proof of Claim \ref{e0}}
Substitute $e=0$ into (\ref{eq1}), then
\begin{eqnarray}
\lambda&=&abd+a(cd+\frac{d^2}{2})+bcd+\frac{2(a+b+c)^3}{25}+\frac{dc^2}{4}. \label{413}
\end{eqnarray}
%subject to $a+b+c+d=1$ and $a, b, d>0$, $c\geq 0.$

If $c=0$, applying (\ref{eqa}), then $b=3a-1$ and $d=2-4a$. By Theorem \ref{KKT}, $\frac{\partial\lambda}{\partial b}=\frac{\partial\lambda}{\partial d}$. Combining with
\begin{eqnarray*}
\frac{\partial\lambda}{\partial b}&=&ad+\frac{6(a+b)^2}{25},\\
\frac{\partial\lambda}{\partial d}&=&ab+ad,
\end{eqnarray*}
we get $ab=\frac{6(a+b)^2}{25}$. Substituting $b=3a-1$, we have $21a^2-23a+6=0$,  then $a=\frac{2}{3}$ or $\frac{3}{7}$. Note that if $a=\frac{2}{3}$, then $d<0$, a contradiction. So $a=\frac{3}{7}$, $b=\frac{2}{7}$ and $d=\frac{2}{7}$. Therefore $\lambda=\frac{4}{49}\leq 0.082.$

If $c\not= 0$, substitute $c=1-a-b-d$ into (\ref{413}),
\begin{eqnarray*}
\lambda&=&abd+a[(1-a-b-d)d+\frac{d^2}{2}]+b(1-a-b-d)d+\frac{2(1-d)^3}{25}+\frac{d(1-a-b-d)^2}{4}\\
&=&ad-a^2d-abd-\frac{ad^2}{2}+bd-b^2d-bd^2+\frac{2(1-d)^3}{25}+\frac{d(1-a-b-d)^2}{4}.
\end{eqnarray*}
subject to $a+b+d\le 1$, and $a. b, d\ge 0$. Since $a+b+d< 1$ and we have showned that $a, b, d>0$, then by Theorem \ref{KKT},
\begin{eqnarray*}
\frac{\partial\lambda}{\partial a}&=&d-2ad-bd-\frac{d^2}{2}-\frac{d(1-a-b-d)}{2}=0,\\
\frac{\partial\lambda}{\partial b}&=&-ad+d-2bd-d^2-\frac{d(1-a-b-d)}{2}=0,\\
\frac{\partial\lambda}{\partial d}&=&a-a^2-ab-ad+b-b^2-2bd-\frac{6(1-d)^2}{25}+\frac{(1-a-b-d)^2}{4}-\frac{d(1-a-b-d)}{2}=0.
\end{eqnarray*}
Then we have $b=1-3a$ (by solving $\frac{\partial\lambda}{\partial a}=0$) and $d=8a-2$ (by solving $\frac{\partial\lambda}{\partial b}=0$), substitute this into $\frac{\partial\lambda}{\partial d}=0$, we have $1266a^2-787a+121=0$. Therefore $a=\frac{787-5\sqrt{265}}{2532}$ and $\lambda\leq0.0939$.
\q

\subsection{\em Proof of Claim \ref{A4}}
Since $G$ contains no other $K_4^3$, then for any $v\in A_4$, there is no $K_3$ in $G_v[A]$, so $G_v[A]$ forms a $C_4$. Hence
$$\lambda(G_v[A], \vec x)\leq max\{(a+b)(c+d), (a+c)(b+d), (a+d)(b+c)\}\leq\frac{(a+b+c+d)^2}{4}.$$
Since $G$ is $K_4^3\cup\{e\}$-free, then $N^*(x, y)\subseteq A$ for $x, y\in [n]\setminus A$. Since $G$ contains only one $K_4^3$, we claim that $|N^*(x, y)|\leq2$ for all $x, y\in A_4$. Recall that $G_x$ and $G_y$ are $C_4$'s, then $G_x[A]\cap G_y[A]$ must be two vertex disjoint edges or $G_x[A]=G_y[A]$. If $|N^*(x, y)|\geq 3$, then there are $z, w\in N^*(x, y)\subseteq A$ such that $zw\in G_x\cap G_y$. Then $\{x, y, z, w\}$ forms a $K_4^3$, a contradiction. Recall that $a\geq b\geq c\geq d$, then we may assume that $N^*(x, y)=\{1, 2\}$ for $x, y\in A_4$ with the Lagrangian non-decreasing. Therefore
\begin{eqnarray*}
\lambda(G[A\cup A_4], \vec x)&\leq&abc+abd+acd+bcd+\frac{(a+b)e^2}{2}+e\frac{(a+b+c+d)^2}{4},
\end{eqnarray*}
subject to $a+b+c+d+e=1$ and $c+d>0.22354$ and $a+b>c+d$. Let $\alpha={a+b\over a+b+c+d+e}$, $\gamma={c+d\over a+b+c+d+e}\geq 0.22354$ and $\eta={e\over a+b+c+d+e}$. So
\begin{eqnarray*}
{\lambda(G[A\cup A_4], \vec x)\over (a+b+c+d+e)^3}&\leq&\frac{\alpha^2\gamma}{4}+\frac{\alpha\gamma^2}{4}+\frac{\alpha\eta^2}{2}+\eta\frac{(\alpha+\gamma)^2}{4}=\lambda(\alpha, \gamma, \eta)\triangleq\lambda
\end{eqnarray*}
subject to
\begin{equation}
\left\{
             \begin{array}{lr}
             \alpha+\gamma+\eta=1, &  \\
             \gamma\geq 0.22354, &\\
             \alpha\geq \gamma. &\label{6}
             \end{array}
\right.
\end{equation}

If $\eta=0$, then $\lambda\leq\lambda(\frac{1}{2}, \frac{1}{2}, 0)=\frac{1}{16}.$

If $\alpha=\gamma$, then $\eta=1-2\alpha$. So $$\lambda=\frac{\alpha^3+\alpha\eta^2}{2}+\eta\alpha^2=\frac{\alpha(\alpha+\eta)^2}{2}=\frac{\alpha(1-\alpha)^2}{2}\leq\frac{(\frac{2\alpha+1-\alpha+1-\alpha}{3})^3}{4}=\frac{2}{27}<0.075.$$

If $\gamma=0.22354$, then $\eta=0.77646-\alpha$
\begin{eqnarray*}
\lambda%&=&\frac{a^3}{4}-\frac{63823a^2}{100000}+\frac{546385151041889892675213a}{1407374883553280000000000}+\frac{4849976047767}{500000000000000}\\
&\leq&\frac{\alpha^3}{4}-0.6382\alpha^2+0.38823\alpha+0.0097=\lambda_0\\
\lambda_0'&=&\frac{3\alpha^2}{4}-1.2764\alpha+0.38823.
\end{eqnarray*}
So $\lambda_0$ is increasing in $[0, 0.3965684 ]$. Then $\lambda\leq 0.0789$.

Therefore we may assume that $\lambda$ gets the maximum in its interior points. By Theorem \ref{KKT}, then $\frac{\partial \lambda}{\partial\gamma}=\frac{\partial \lambda}{\partial\eta}$. Combining with
\begin{eqnarray*}
%\frac{\partial \lambda}{\partial a}&=&\frac{ac}{2}+\frac{c^2}{4}+\frac{e^2}{2}+e\frac{(a+c)}{2}\\
\frac{\partial \lambda}{\partial\gamma}&=&\frac{\alpha^2}{4}+\frac{\alpha\gamma}{2}+\eta\frac{(\alpha+\gamma)}{2}\\
\frac{\partial \lambda}{\partial\eta}&=&\alpha\eta+\frac{(\alpha+\gamma)^2}{4},
\end{eqnarray*}
we get $\gamma^2=2\eta(\gamma-\alpha)$, contradicting to $\alpha>\gamma$. This completes the proof of Claim \ref{A4}.
\q

\subsection{\em Proof of Claim \ref{A3}.}
Note that for $x, y\in A_3$, if $|N^*(x, y)|=4$ (i.e. $N^*(x, y)=A$) and there exists $zw\in G_x[A]\cap G_y[A]$, then $\{x, y, z, w\}$ forms a $K_4^3$, a contradiction. So for $x, y\in A_3$, if $|N^*(x, y)|=4$, then $G_x[A]\cap G_y[A]=\emptyset$ and $G_x[A]\cup G_y[A]={A\choose 2}$. Such pairs $\{x, y\}$ are partitioned into at most ${6\choose 3}/2$ groups $\{B_{0i}, B_{1i}\}$ such that $G_x[A]$ are the same for all $x\in B_{0i}$, $G_y[A]$ are the same for all $y\in B_{1i}$, and $G_x[A]\cup G_y[A]$ is a partition of ${A\choose 2}$ for all $x\in B_{0i}$ and $y\in B_{1i}$ for each $i$. Assume that there are $s$ such groups.
So $|N^*(x, y)|\leq 3$ for $x, y\in A_4\cup A_3$ except $x\in B_{0i}$ and $y\in B_{1i}$, then we may assume that $N^*(x, y)=\{1, 2, 3\}$ for such $x, y$ with the Lagrangian non-decreasing. Note that $|N^*(x, y)|\leq 2$ for $x, y\in B_{0i}$ or $x, y\in B_{1i}$ for $1\leq i\leq s$. Therefore we may assume that $N^*(x, y)=\{1, 2\}$ for such $x, y$ with the Lagrangian non-decreasing.
And for all $x\in A_3$, if $x23\notin E(G)$, then we may assume that $x12, x13, x14\in E(G)$. If $x23\in E(G)$, since $\{x, 1, 2, 3\}$ doesn't span $K_4^3$, then one of $x12, x13$ does not belong to $E(G)$, we can replace $x23$ by that edge and replace other 2 edges $xij$, $ij\in A^{(2)}$ from $\{x12, x13, x14\}$ with the Lagrangian non-decreasing. Let $f_{0i}=\sum_{v\in B_{0i}}x_v$ and $f_{1i}=\sum_{v\in B_{1i}}x_v$ and $f_i=f_{0i}+f_{1i}$ for $1\leq i\leq s$. Let $f'=f-\sum_{i=1}^sf_i$. By Claim \ref{A4} and the above analysis, we have
\begin{eqnarray*}
\lambda(G[A\cup A_4\cup A_3], \vec x)&\leq&0.0789(a+b+c+d+e)^3+(ab+ac+ad)f+(ef+f'\sum_{i=1}^sf_i+\sum_{1\leq i\not=j\leq s}f_if_j\\
&+&\frac{f'^2}{2})(a+b+c)+\sum_{i=1}^s\frac{f_{0i}^2+f_{1i}^2}{2}(a+b)+\sum_{i=1}^sf_{0i}f_{1i}(a+b+c+d)\\
&\leq&0.0789(a+b+c+d+e)^3+(ab+ac+ad)f+ef(a+b+c)+{(a+b)f^2\over 2}\\
&+&c(f'\sum_{i=1}^sf_i+\sum_{1\leq i\not=j\leq s}f_if_j+\frac{f'^2}{2})+(c+d)\sum_{i=1}^sf_{0i}f_{1i}\\
&\leq&0.0789(a+b+c+d+e)^3+(ab+ac+ad)f+ef(a+b+c)+{(a+b)f^2\over 2}\\
&+&c(f'\sum_{i=1}^sf_i+\sum_{1\leq i\not=j\leq s}f_if_j+\frac{f'^2}{2})+(c+d)\sum_{i=1}^s{f_i^2\over 4}\\
&=&\lambda(f',f_1,\dots,f_s).
\end{eqnarray*}
Note that
\begin{eqnarray*}
\lambda(f, 0, \dots, 0)&-&\lambda(f',f_1,\dots,f_s)\\
&=&\frac{cf^2}{2}-c(f'\sum_{i=1}^sf_i+\sum_{1\leq i\not=j\leq s}f_if_j+\frac{f'^2}{2})-(c+d)\sum_{i=1}^s{f_i^2\over 4}\\
&=&c{\sum_{i=1}^sf_i^2\over 2}-(c+d)\sum_{i=1}^s{f_i^2\over 4}\\
&\geq&0
\end{eqnarray*}
since $c\geq d$. So we may assume that $f'=f$ and $f_i=0$ for $1\leq i\leq s$. So
\begin{eqnarray}
\lambda(G[A\cup A_4\cup A_3], \vec x)&\leq&0.0789(a+b+c+d+e)^3+(ab+ac+ad)f+(ef+\frac{f^2}{2})(a+b+c)\notag \\
&=&0.0789(a+(b+c)+d+e)^3+(a(b+c)+ad)f+(ef+\frac{f^2}{2})(a+(b+c)) \notag \\
&=&\lambda.
\end{eqnarray}
Let $\alpha={a\over (a+b+c+d+e+f)}$, $\beta={b+c\over a+b+c+d+e+f}$, $\delta={d\over a+b+c+d+e+f}$, $\eta={e\over a+b+c+d+e+f}$ and $\rho={f\over a+b+c+d+e+f}$, and let $\tau={\lambda\over (a+b+c+d+e+f)^3}$. Then
$$\tau=0.0789(\alpha+\beta+\delta+\eta)^3+(\alpha\beta+\alpha\delta)\rho+(\eta\rho+\frac{\rho^2}{2})(\alpha+\beta),$$
and it's sufficient to prove that $\tau\leq 0.092.$

{\em Case 1.} $d\geq0.11177$.

In this case, note that $\delta\geq d\geq 0.11177$. Recall that $b+c\geq 2d$, then $\beta\geq b+c\geq0.22354$. Note that $\tau$ is non-decreasing if we change $(\beta, \delta)$ to $(\beta+\epsilon, \delta-\epsilon)$ for $\epsilon>0$. So, we may assume that $\delta=0.11177$. By Claim \ref{2}, $a\leq {3-\sqrt 3\over 3}-d\leq 3\times 0.11177\leq b+c+d$, so $\alpha\leq\beta+\delta$. Therefore $\tau$ is non-decreasing if we change $(\alpha, \beta)$ to $(\alpha+\epsilon, \beta-\epsilon)$ for small $\epsilon$, we may assume that $\beta=0.22354$. Therefore
\begin{eqnarray}
\tau&=&0.0789(\alpha+0.33531+\eta)^3+0.33531\alpha\rho+(\eta\rho+\frac{\rho^2}{2})(\alpha+0.22354)
\end{eqnarray}
subject to $\alpha+0.33531+\eta+\rho=1$, $\alpha, \eta, \rho \geq 0$. If $\rho=0,$ then $\tau=0.0789$. So we may assume that $\rho>0$.

If $\alpha=0$, then $\eta=0.66469-\rho$. So
\begin{eqnarray*}
\tau&\leq&-0.0789\rho^3+0.125\rho^2-0.08\rho+0.0789=f(\rho),\\
f'(\rho)&=&-0.2367\rho^2+0.25\rho-0.08<0.
\end{eqnarray*}
Note that $f(\rho)$ is decreasing in $[0, 1]$, then $\tau_0\leq f(0)\leq 0.0789.$ So we may assume that $\alpha>0$.

If $\eta>0$, then by Theorem \ref{KKT}, we have $\frac{\partial \tau}{\partial \alpha}=\frac{\partial \tau}{\partial \eta}$, so $\alpha=0.11177+\eta+\frac{\rho}{2}$. Therefore $\rho=1.55292-4\alpha$ and $\eta=3\alpha-0.88823$. So
\begin{eqnarray*}
\tau%&=&0.0789(a+b+d+e)^3+(ab+ad)f+(ef+\frac{f^2}{2})(a+b)\\
%&=&\frac{656a^3}{625}-\frac{224956541a^2}{62500000}+\frac{38529999873743a}{25000000000000}-\frac{940906646068298163}{10000000000000000000}\\
&\leq&1.05\alpha^3-3.55\alpha^2+1.55\alpha-0.094=\tau_0\\
\tau_0'&=&3.15\alpha^2-7.1\alpha+1.55.
\end{eqnarray*}
Note that $\tau_0$ is increasing in $[0, \frac{71-4\sqrt {193}}{63}]$, then $\tau_0\leq 0.09.$

If $\eta=0$, then $\rho=0.66469-\alpha>0$. So
\begin{eqnarray*}
\tau&=&0.0789(\alpha+0.33531)^3+0.33531\alpha(0.66469-\alpha)+\frac{(0.66469-\alpha)^2}{2}(\alpha+0.22354)\\
&\leq&0.5789\alpha^3-0.8088\alpha^2+0.3219\alpha+0.0524=\tau_0\\
\tau_0'&=&1.7367\alpha^2-1.6176\alpha+0.3219.
\end{eqnarray*}
Note that $\tau_0$ is increasing in $[0, \frac{2696-\sqrt {1056819}}{5789}]$ or $[\frac{2696+\sqrt {1056819}}{5789}, 0.66469]$, then $\tau_0\leq 0.092.$

{\em Case 2.} $d<0.11177.$

In this case, we know that $b+c\geq 0.307$, so $\beta\geq 0.307$. By Claim \ref{3}, we know that $d\geq0.0848$, so $\delta\geq 0.0848$. Note that $\tau$ is non-decreasing if we change $(\beta, \delta)$ to $(\beta+\epsilon, \delta-\epsilon)$ for $\epsilon>0$, so we may let $\delta=0.0848$ in $\tau$. Therefore
\begin{eqnarray*}
\tau&=&0.0789(\alpha+\beta+0.0848+\eta)^3+(\alpha\beta+0.0848\alpha)\rho+(\eta\rho+\frac{\rho^2}{2})(\alpha+\beta)
\end{eqnarray*}
subject to
\begin{equation}
\left\{
             \begin{array}{lr}
             \alpha+\beta+0.0848+\eta+\rho=1, &  \\
             \beta\geq0.307, &\label{9} \\
             \alpha, \rho\geq0 &
             \end{array}
\right.
\end{equation}
Note that $\rho>0$.

If $\alpha=0$, then we claim that $\beta=0.307$, this is because that $\tau$ is non-decreasing if we change $(\alpha, \beta)$ to $(\alpha+\epsilon, \beta-\epsilon)$ for small $\epsilon$. So $\eta=0.6082-\rho$. Substituting these into (\ref{9}), we have
\begin{eqnarray*}
\tau&\leq&-0.0789\rho^3+0.042\rho^2+0.002\rho+0.0789=f(\rho),\\
f'(\rho)&=&-0.2367\rho^2+0.084\rho+0.002.
\end{eqnarray*}
Note that $f(\rho)$ is increasing in $[0, {2\sqrt{6215}+140\over 789}]$ and decreasing in $[{2\sqrt{6215}+140\over 789}, 1]$, then $\tau\leq f({2\sqrt{6215}+140\over 789})< 0.085.$ So we may assume that $\alpha>0$.

If $\beta>0.307$, then by Theorem \ref{KKT}, $\frac{\partial\tau}{\partial\alpha}=\frac{\partial\tau}{\partial\beta}$, simplifying it,  we have $\alpha=\beta+0.0848>0.3918$. If $\eta>0$, then  by Theorem \ref{KKT}, we have $\frac{\partial\tau}{\partial\alpha}=\frac{\partial\tau}{\partial\eta}$, simplyfying it, we get $\alpha=0.0848+\eta+\frac{\rho}{2}.$  Hence $1.1754<3\alpha=\alpha+\beta+2\times 0.0848+\eta+\frac{\rho}{2}<1+0.0848=1.0848$, a contradiction. So $\eta=0$. Then $\rho=1-\alpha-\beta-0.0848=1-2\alpha>0$ and $\beta=\alpha-0.0848.$ So $0.39<\alpha<0.5$. And
\begin{eqnarray*}
\tau&=&0.0789(2\alpha)^3+\alpha^2(1-2\alpha)+\frac{(1-2\alpha)^2}{2}(2\alpha-0.0848)\\
&\leq&2.65\alpha^3-3.1695\alpha^2+1.1695\alpha-0.0424=\tau_0,\\
\tau_0'&=&7.95\alpha^2-6.339\alpha+1.1695.
\end{eqnarray*}
Note that $\tau_0$ is decreasing in $[0.39, 0.5]$, then $\tau_0\leq \tau_0(0.39)\le 0.089.$

So we may assume that $\beta=0.307$. If $\eta>0$, then by Theorem \ref{KKT}, $\frac{\partial\tau}{\partial\alpha}=\frac{\partial\tau}{\partial\eta}$, simplifying it, we have $\alpha=0.0848+\eta+\frac{\rho}{2}.$ Since $\alpha+\eta+\rho=1-\beta-0.0848=0.6082$, then $\eta=3\alpha-0.7778$ and $\rho=1.386-4\alpha>0.$ Since $\frac{\partial\tau}{\partial\eta}=\frac{\partial\tau}{\partial\rho}$, then $0.2367(\alpha+\beta+0.0848+\eta)^2=\alpha\beta+0.0848\alpha+\eta(\alpha+\beta)$. By direct calculation, $\frac{492\alpha^2}{625}-\frac{395603\alpha}{312500}+\frac{685129883}{2500000000}=0$,  and $\alpha=\frac{395603-20\sqrt{180576895}}{492000}$. But $3\alpha<0.774$, then $\eta=3\alpha-0.7778<0$, a contradiction. So $\eta=0$, then $\rho=1-\alpha-\beta-0.0848=0.6082-\alpha>0.$ So
\begin{eqnarray*}
\tau&=&0.0789(\alpha+0.3918)^3+0.3918\alpha(0.6082-\alpha)+\frac{(0.6082-\alpha)^2}{2}(\alpha+0.307)\\
&\leq&0.5789\alpha^3-0.75375\alpha^2+0.27286\alpha+0.06153=\tau_0,\\
\tau_0'&=&1.7367\alpha^2-1.5075\alpha+0.27286.
\end{eqnarray*}
Note that $\tau_0$ is increasing in $[0, \frac{5025}{11578}-\frac{\sqrt{942631005}}{173670}]$, then $\tau_0\leq 0.092.$ This complete the proof of Claim \ref{A3}.
\q

\subsection{\em Proof of Claim \ref{A2}}
By Claim \ref{A3}, we have
\begin{eqnarray*}
\lambda(G)&\leq&0.092(a+b+c+d+e+f)^3+(ab+ac)g+abh+(eg+eh+fg+fh+\frac{g^2}{2}+gh\\
&+&\frac{h^2}{2})(a+b+c+d)\\
&\leq&0.092(a+\beta+d+\zeta)^3+a\beta\eta+(\zeta\eta+\frac{\eta^2}{2})(a+\beta+d)\\
&=&\lambda,
\end{eqnarray*}
where $\eta=g+h$, $\beta=b+c$ and $\zeta=e+f$.

{\em Case 1.} $d\geq 0.11177.$

Since replacing $(a, d)$ by $(a+\epsilon, d-\epsilon)$ for $\epsilon>0$ will not decrease $\lambda$, so we may assume that $d=0.11177$. Let $\alpha=a+\beta$. Then
\begin{eqnarray*}
\lambda&\leq&0.092(\alpha+\delta+\zeta)^3+\frac{\alpha^2\eta}{4}+(\zeta\eta+\frac{\eta^2}{2})(\alpha+\delta),
\end{eqnarray*}
subject to
\begin{equation}
\left\{
             \begin{array}{lr}
             \alpha+\delta+\zeta+\eta=1, &  \\
             \delta=0.11177, \alpha, \zeta, \eta\ge 0. &\label{10}
             \end{array}
\right.
\end{equation}

If $\eta=0$, then $\lambda=0.092$, we are done. So assume that $\eta>0.$

If $\alpha=0$, then $\zeta=0.88823-\eta$.
\begin{eqnarray*}
\lambda&\leq&-0.092\eta^3+ 0.221\eta^2-0.1767\eta+0.092=\lambda_0\\
\lambda_0'&=&-0.276\eta^2+0.442\eta-0.1767.
\end{eqnarray*}
Note that $\lambda_0$ is decreasing in $[0, 0.7700236]$ or $[0.8314257, 1]$. Therefore $\lambda_0\leq 0.095.$ So assume that $\alpha>0.$

If $\zeta=0$, then $\eta=0.88823-\alpha>0$. So
\begin{eqnarray*}
\lambda%&=&0.092(a+0.11177)^3+\frac{a^2(0.88823-a)}{4}+\frac{(0.88823-a)^2}{2}(a+0.11177)\\
%&=&\frac{171a^3}{500}-\frac{28971949a^2}{50000000}+\frac{186654211519a}{625000000000}+\frac{2763691938632621}{62500000000000000}\\
&\leq&0.342\alpha^3-0.58\alpha^2+0.3\alpha+0.045=\lambda_0\\
\lambda_0'&=&1.026\alpha^2-1.16\alpha+0.3.
\end{eqnarray*}
Note that $\lambda_0$ is increasing in $[0, \frac{290-5\sqrt {286}}{513}]$ or $[\frac{290+5\sqrt {286}}{513}, 0.88823]$. Therefore $\lambda_0\leq 0.095.$

If $\zeta>0$, then by Theorem \ref{KKT}, we have $\frac{\partial\lambda}{\partial \alpha}=\frac{\partial\lambda}{\partial \zeta}$, i.e. $2\zeta+\eta=\alpha+2\delta$. Since $\alpha+\delta+\zeta+\eta=1$, then $\zeta=2\alpha-0.66469$ and $\eta=1.55292-3\alpha>0$. So $0.33234<\alpha<0.51764$. And
\begin{eqnarray*}
\lambda%&=&0.092(3a-0.55292)^3+\frac{a^2(1.55292-3a)}{4}+((2a-0.66469)(1.55292-3a)\\
%&+&\frac{(1.55292-3a)^2}{2})(a+0.11177)\\
%&=&\frac{117a^3}{500}-\frac{17793207a^2}{25000000}+\frac{595017249699a}{1250000000000}+\frac{120259899781747}{31250000000000000}\\
&\leq&0.234\alpha^3-0.7115\alpha^2+0.4765\alpha+0.00385=\lambda_0\\
\lambda_0'&=&0.702\alpha^2-1.423\alpha+0.4765.
\end{eqnarray*}
Note that $\lambda_0$ is increasing in $[0, \frac{1423-11\sqrt{5677}}{1404}]$. Therefore $\lambda_0\leq 0.096.$

{\em Case 2.} $d<0.11177.$

In this case, we have $\beta\geq0.307$. By Claim \ref{3}, then $d\geq0.0848$. Since replacing $(a, d)$ by $(a+\epsilon, d-\epsilon)$ for $\epsilon>0$ will not decrease $\lambda$, so we may assume $d=0.0848$. So
\begin{eqnarray*}
\lambda&=&0.092(a+\beta+\delta+\zeta)^3+a\beta\eta+(\zeta\eta+\frac{\eta^2}{2})(a+\beta+\delta),
\end{eqnarray*}
subject to
\begin{equation}
\left\{
             \begin{array}{lr}
             a+\beta+\delta+\zeta+\eta=1, &  \\
             \beta\geq0.307, &   \\
             \delta=0.0848. &\label{11}
             \end{array}
\right.
\end{equation}

If $a=0$, we claim that $\beta=0.307$ since $\lambda$ is non-decreasing if we change $(a, \beta)$ to $(a+\epsilon, \beta-\epsilon)$ for small $\epsilon$. So $\zeta=0.6082-\eta$. Then
\begin{eqnarray*}
\lambda&\leq&-0.092\eta^3+0.0801\eta^2-0.037\eta+0.092=f(\eta),\\
f'(\eta)&=&-0.276\eta^2+0.1602\eta-0.037<0.
\end{eqnarray*}
Therefore $\lambda_0\leq f(0)=0.092$. So we may assume that $a>0$.

If $\beta>0.307$, then by Theorem \ref{KKT}, $\frac{\partial\lambda}{\partial a}=\frac{\partial\lambda}{\partial \beta}$, simplifying it, we get $a=\beta>0.307$. If $\zeta>0$, then we have $\frac{\partial\lambda}{\partial a}=\frac{\partial\lambda}{\partial \zeta}$, i.e. $a+\delta=\zeta+\frac{\eta}{2}.$ So $1<3a+2\delta=a+\beta+\delta+\zeta+\frac{\eta}{2}<1$, a contradiction. Then $\zeta=0$. So $\eta=1-a-\beta-\delta=1-2a-\delta=0.9152-2a>0$, then $a<0.4576$. Recall that $\delta=0.0848$. So
\begin{eqnarray*}
\lambda&=&0.092(2a+0.0848)^3+(0.9152-2a)a^2+\frac{(0.9152-2a)^2}{2}(2a+0.0848).%\\
%&=&\frac{342a^3}{125}-\frac{193936a^2}{78125}+\frac{33512821a}{48828125}+\frac{2171018171}{61035156250}=\lambda_0\\
%\lambda_0'&=&\frac{1026a^2}{125}-\frac{387872a}{78125}+\frac{33512821}{48828125}.
\end{eqnarray*}
By a direct calculation on the derivative of $\lambda_0$, we obtain that $\lambda_0$ is increasing in $[0, 0.2138466]$. Therefore $\lambda_0\leq 0.096.$

So we may assume that $\beta=0.307$ and $\delta=0.0848$. If $\zeta=0$, then $\eta=1-a-\beta-\delta=0.6082-a$. By Theorem \ref{KKT}, then $\frac{\partial\lambda}{\partial a}=\frac{\partial\lambda}{\partial \eta}$. Combining them, we get
$0.276(a+0.3918)^2+\frac{(0.6082-a)^2}{2}=(0.6082-a)(a+0.0848)+0.307a$. Solving the equation for $a$ and substituting the values  into $\lambda$, we obtain that $\lambda\leq 0.095.$ So we may assume that $\zeta>0$. By Theorem \ref{KKT}, then $\frac{\partial\lambda}{\partial a}=\frac{\partial\lambda}{\partial \zeta}$, i.e. $a+\delta=\zeta+\frac{\eta}{2}$. So $\eta=1.0468-4a$ and $\zeta=3a-0.4386.$ Then $0.1462<a<0.2617$ and
\begin{eqnarray*}
\lambda&=&0.092(4a-0.0468)^3+0.307(1.0468-4a)a+((3a-0.4386)(1.0468-4a)\\
&+&\frac{(1.0468-4a)^2}{2})(a+0.3918).%\\
%&=&\frac{236a^3}{125}-\frac{716959a^2}{312500}+\frac{538899957a}{781250000}+\frac{135820792401}{3906250000000},\\
%\lambda'&=&\frac{708a^2}{125}-\frac{716959a}{156250}+\frac{538899957}{781250000}.
\end{eqnarray*}
By a direct calculation on the derivative of $\lambda$, we obtain that $\lambda$ is increasing in $[0, \frac{716959}{1770000}-\frac{\sqrt{211982461}}{70800}]$ or $[\frac{716959}{1770000}+\frac{\sqrt{211982461}}{70800}, 1]$. Recall that $0.1462<a<0.2617$, then  $\lambda<0.0961$. \q

\subsection{\em Proof of Claim \ref{2f0}}
If $f=0$, then
\begin{eqnarray*}
\lambda&=&\frac{\alpha^2(c+d+e)}{4}+\alpha c(d+e)+\alpha(de+\frac{e^2}{2})+\frac{d^2(\alpha+e)}{4}
\end{eqnarray*}
subject to
\begin{equation}
\left\{
             \begin{array}{lr}
             \alpha+c+d+e=1, &  \\
             c\geq0.08. &
             \end{array}
\right.
\end{equation}

If $d=0$, then $\lambda=\frac{\alpha^2(c+e)}{4}+\alpha ce+\alpha\frac{e^2}{2}\leq \lim_{n\to\infty}B(2, n-2)=\frac{\sqrt 3}{18}$. So we are done.

If $e=0$, then $\lambda=\frac{\alpha^2(c+d)}{4}+\alpha cd+\frac{\alpha d^2}{4}\leq \lambda(K_5^3)<\frac{\sqrt 3}{18}$. So we are done.

So we may assume that $d, e>0$. By direct calculation,
\begin{eqnarray*}
\frac{\partial \lambda}{\partial \alpha}&=&\frac{\alpha(c+d+e)}{2}+cd+ce+de+\frac{e^2}{2}+\frac{d^2}{4}\\
\frac{\partial \lambda}{\partial c}&=&\frac{\alpha^2}{4}+\alpha d+\alpha e\\
\frac{\partial \lambda}{\partial d}&=&\frac{\alpha^2}{4}+\alpha c+\alpha e+\frac{d(\alpha+e)}{2}\\
\frac{\partial \lambda}{\partial e}&=&\frac{\alpha^2}{4}+\alpha c+\alpha d+\alpha e+\frac{d^2}{4}.
\end{eqnarray*}
By Theorem \ref{KKT},  $\frac{\partial \lambda}{\partial d}=\frac{\partial \lambda}{\partial e}$, combining with the above equations, we get $d=2e-2\alpha.$ If $c>0.08$, then by Theorem \ref{KKT},  $\frac{\partial \lambda}{\partial c}=\frac{\partial \lambda}{\partial e}$, combining the equations, we get $\alpha c+{d^2\over 4}=0$, a contradiction. So $c=0.08$,  therefore $\alpha+d+e=0.92$, combining with $d=2e-2\alpha$, we get $d=\frac{1.84-4\alpha}{3}$ and $e=\frac{\alpha+0.92}{3}$. Since $\frac{\partial \lambda}{\partial \alpha}=\frac{\partial \lambda}{\partial e}$, then $\frac{5\alpha^2}{36}+\frac{307\alpha}{450}-\frac{3473}{11250}=0$. Then $\alpha=\frac{3\sqrt{14331}-307}{125}$. By direct calculation, $\lambda\le 0.096<\frac{\sqrt 3}{18}$, a contradiction. \q

\subsection{\em Proof of Claim \ref{2e0}}
If $e=0$, in view of (\ref{eq6}), then
\begin{eqnarray*}
\lambda&=&\frac{\alpha^2(c+d+f)}{4}+\alpha cd+(\alpha+c)(df+\frac{f^2}{2})+\frac{d^2(\alpha+f)}{4},
\end{eqnarray*}
subject to
\begin{equation}
\left\{
             \begin{array}{lr}
             \alpha+c+d+f=1, &  \\
             c\geq0.08; \alpha, d, f \ge 0. &
             \end{array}
\right.
\end{equation}

If $d=0$, then we have
\begin{eqnarray*}
\lambda&=&\frac{\alpha^2(c+f)}{4}+(\alpha+c)\frac{f^2}{2}.
\end{eqnarray*}
We relax the constraint $c\geq 0.08$ to $c\geq 0$ and $\alpha+c+f=1$. If $c=0$, then $\lambda=\frac{\alpha^2(1-\alpha)}{4}+\alpha\frac{(1-\alpha)^2}{2}\leq\frac{\sqrt 3}{18}$ by Fact \ref{aaa}. If $c>0$, recall that $f>0$, then by Theorem \ref{KKT}, we have $\frac{\partial \lambda}{\partial c}=\frac{\partial \lambda}{\partial f}$, then $f=2(\alpha+c)$, combining with  $\alpha+c+f=1$, we have $f=\frac{2}{3}$ and $c=\frac{1}{3}-\alpha$. So $\lambda=\frac{\alpha^2(1-\alpha)}{4}+\frac{2}{27}$, where $\alpha<\frac{1}{3}$. Since $\lambda'=\frac{\alpha(2-3\alpha)}{4}$, then $\lambda$ is increasing in $[0, \frac{1}{3}]$, i.e. $\lambda\leq\frac{(\frac{1}{3})^2\frac{2}{3}}{4}+\frac{2}{27}=\frac{5}{54}\leq\frac{\sqrt 3}{18}.$

So we may assume that $d>0$. Recall that $e+f>\alpha+\frac{d}{2}$ and
$e=0$, so $f>\alpha+\frac{d}{2}$. By direct calculation,
\begin{eqnarray*}
\frac{\partial \lambda}{\partial \alpha}&=&\frac{\alpha(c+d+f)}{2}+cd+df+\frac{f^2}{2}+\frac{d^2}{4}\\
\frac{\partial \lambda}{\partial c}&=&\frac{\alpha^2}{4}+\alpha d+df+\frac{f^2}{2}\\
\frac{\partial \lambda}{\partial d}&=&\frac{\alpha^2}{4}+\alpha c+(\alpha+c)f+\frac{d(\alpha+f)}{2}\\
\frac{\partial \lambda}{\partial f}&=&\frac{\alpha^2}{4}+(\alpha+c)(d+f)+\frac{d^2}{4}.
\end{eqnarray*}
By Theorem \ref{KKT}, we have $\frac{\partial \lambda}{\partial d}=\frac{\partial \lambda}{\partial f}$, i.e. $$\frac{d^2}{4}+(\alpha+c)d=\alpha c+\frac{d(\alpha+f)}{2}>\alpha c+\frac{d(2\alpha+\frac{d}{2})}{2},$$ then $d>\alpha.$
If $c>0.08$, then $\frac{\partial \lambda}{\partial \alpha}=\frac{\partial \lambda}{\partial c}$.
Therefore $\frac{\alpha^2}{4}+\alpha d=\frac{\alpha(c+d+f)}{2}+cd+\frac{d^2}{4}$. Since $f>\alpha+\frac{d}{2},$ then $\frac{\alpha^2}{4}+\alpha d>\frac{\alpha(c+\alpha+\frac{3d}{2})}{2}+cd+\frac{d^2}{4}>\frac{\alpha^2}{2}+\frac{d^2+3\alpha d}{4}$. So $\alpha d>\alpha^2+d^2$, a contradiction.
So $c=0.08$. Since $\frac{\partial \lambda}{\partial \alpha}=\frac{\partial \lambda}{\partial f}$ and $f=0.92-\alpha-d$, then $$\frac{\alpha^2}{4}+\alpha d+\alpha f+cf-\frac{\alpha(1-\alpha)}{2}-df-\frac{f^2}{2}=0.$$
Substituting $c=0.08$ and $f=0.92-\alpha-d$ into it, we have
$$\frac{\alpha^2}{4}+\alpha d+\alpha(0.92-\alpha-d)+0.08(0.92-\alpha-d)-\frac{\alpha(1-\alpha)}{2}-d(0.92-\alpha-d)-\frac{(0.92-\alpha-d)^2}{2}=0.$$
Simplifying it, we have $\frac{d^2}{2}-\frac{2d}{25}+\frac{63\alpha}{50}-\frac{3\alpha^2}{4}-\frac{437}{1250}=0$, then $d=\sqrt{\frac{3\alpha^2}{2}-\frac{63\alpha}{25}+\frac{441}{625}}+\frac{2}{25}$. Note that $0.92=\alpha+d+f>\alpha+\alpha+\frac{3\alpha}{2}=\frac{7\alpha}{2}$, then $0<\alpha<0.2629$. Therefore $\phi(\alpha)=\frac{3\alpha^2}{2}-\frac{63\alpha}{25}+\frac{441}{625}>\phi(0.2629)>0.1467$. So $d>\sqrt{0.1467}+\frac{2}{25}>0.46$ and $f>0.23+\alpha$, then $\alpha+d+f>\alpha+\sqrt{\frac{3\alpha^2}{2}-\frac{63\alpha}{25}+\frac{441}{625}}+\frac{2}{25}+0.23+\alpha\triangleq\psi(\alpha)>\psi(0)=1.15>1$, a contradiction. \q

\subsection{\em Proof of Claim \ref{2d0}}
If $d=0$, then we have $\alpha+c+e+f=1$ and $\alpha=e$. In view of (\ref{eq6}),
\begin{eqnarray*}
\lambda&=&\frac{\alpha^2(c+e+f)}{4}+\alpha ce+\alpha\frac{e^2}{2}+(\alpha+c)(ef+\frac{f^2}{2}).
\end{eqnarray*}
Then
\begin{eqnarray*}
\frac{\partial \lambda}{\partial \alpha}&=&\frac{\alpha(c+e+f)}{2}+ce+\frac{e^2}{2}+ef+\frac{f^2}{2}\\
\frac{\partial \lambda}{\partial c}&=&\frac{\alpha^2}{4}+\alpha e+ef+\frac{f^2}{2}\\
\frac{\partial \lambda}{\partial e}&=&\frac{\alpha^2}{4}+\alpha c+\alpha e+(\alpha+c)f,
\end{eqnarray*}

If $c>0.08$, then by Theorem \ref{KKT}, we have $\frac{\partial\lambda}{\partial \alpha}=\frac{\partial\lambda}{\partial c}$, so $\frac{\alpha(1-\alpha)}{2}+\alpha c+\frac{\alpha^2}{2}=\frac{\alpha^2}{4}+\alpha^2$, then $0.08<c=\frac{5\alpha-2}{4}$ and $0<f=\frac{6-13\alpha}{4}$. So $0.464<\alpha<\frac{6}{13}<0.462$, a contradiction.
%Since $\frac{\partial\lambda}{\partial c}=\frac{\partial\lambda}{\partial e}$, then $\frac{f^2}{2}=ac+cf$, i.e. $259a^2-252a+60=0$. So $a=\frac{126\pm4\sqrt{21}}{259}$. Since $f>0$, then $a<\frac{6}{13}$. So $a=\frac{126-4\sqrt{21}}{259}$. However $c<0.02$, a contradiction.

If $c=0.08$, then $f=0.92-2\alpha$. So
\begin{eqnarray*}
\lambda&=&\frac{\alpha^2(1-\alpha)}{4}+0.08\alpha^2+\frac{\alpha^3}{2}+(\alpha+0.08)(\alpha(0.92-2\alpha)+\frac{(0.92-2\alpha)^2}{2})\\
&=&\frac{\alpha^3}{4}-\frac{59\alpha^2}{100}+\frac{437\alpha}{1250}+\frac{529}{15625}.\\
\lambda'&=&\frac{3\alpha^2}{4}-\frac{59\alpha}{50}+\frac{437}{1250}.
\end{eqnarray*}
Note that $\lambda$ is increasing in $[0, \frac{59-\sqrt{859}}{75}]$, then $\lambda\leq0.0955$.
\q

\bigskip

{\bf Acknowledgement}
We are grateful to two reviewers for checking all the details and  giving us valuable  comments to help improve the presentation. The research is supported in part by National Natural Science Foundation of China (No. 11931002).

\end{document}